\documentclass[trsc]{informs3}

\OneAndAHalfSpacedXI 

\usepackage{natbib}
 \bibpunct[, ]{(}{)}{,}{a}{}{,}%
 \def\bibfont{\small}%
\TheoremsNumberedThrough     

\EquationsNumberedThrough    

\MANUSCRIPTNO{TS-2023-0111.R1} 

\usepackage{booktabs}
\usepackage{multirow}
\usepackage{subcaption}
\usepackage{graphicx}
\usepackage{xcolor}
\usepackage{tikz}

\usepackage{mathtools}

\usepackage{algorithm}
\usepackage{algpseudocode}
\usepackage{subfiles} 
\usepackage{bbm}
\usepackage{bookmark}
\usepackage{xr}

\newcommand{\dynamicproblem}{DDWP}
\newcommand{\sigmap}{\sigma^{\prime}}
\newcommand{\sdepot}{0}
\newcommand{\edepot}{|\mathcal{N}| + 1}

\begin{document}


\RUNAUTHOR{Lan et al.}

\RUNTITLE{An iterative sample scenario approach for the DDWP}

\TITLE{An iterative sample scenario approach for the dynamic dispatch waves problem}

\ARTICLEAUTHORS{%
\AUTHOR{Leon Lan, Jasper van Doorn}
\AFF{Department of Mathematics, Vrije Universiteit Amsterdam, \EMAIL{\{l.lan, j.m.h.van.doorn\}@vu.nl}, \URL{}}
\AUTHOR{Niels A. Wouda, Arpan Rijal}
\AFF{Department of Operations, University of Groningen, \EMAIL{\{n.a.wouda, a.rijal\}@rug.nl}, \URL{}}
\AUTHOR{Sandjai Bhulai}
\AFF{Department of Mathematics, Vrije Universiteit Amsterdam, \EMAIL{s.bhulai@vu.nl}. \URL{}}
} 

\ABSTRACT{%
A challenge in same-day delivery operations is that delivery requests are typically not known beforehand, but are instead revealed dynamically during the day.
This uncertainty introduces a trade-off between dispatching vehicles to serve requests as soon as they are revealed to ensure timely delivery, and delaying the dispatching decision to consolidate routing decisions with future, currently unknown requests.
In this paper, we study the \emph{dynamic dispatch waves} problem, a same-day delivery problem in which vehicles are dispatched at fixed decision moments.
At each decision moment, the system operator must decide which of the known requests to dispatch, and how to route these dispatched requests.
The operator's goal is to minimize the total routing cost while ensuring that all requests are served on time.
We propose \emph{iterative conditional dispatch} (ICD), an iterative solution construction procedure based on a sample scenario approach.
ICD iteratively solves sample scenarios to classify requests to be dispatched, postponed, or undecided.
The set of undecided requests shrinks in each iteration until a final dispatching decision is made in the last iteration.
We develop two variants of ICD: one variant based on thresholds, and another variant based on similarity.
A significant strength of ICD is that it is conceptually simple and easy to implement.
This simplicity does not harm performance: through rigorous numerical experiments, we show that both variants efficiently navigate the large state and action spaces of the dynamic dispatch waves problem and quickly converge to a high-quality solution.
Finally, we demonstrate that ICD achieves excellent results on instances from the \emph{EURO meets NeurIPS 2022} vehicle routing competition, nearly matching the performance of the winning machine learning-based strategy.

}%


\maketitle


\section{Introduction}
E-commerce retailers are increasingly offering same-day or even instant delivery services, with deliveries often occurring within two hours after order placement~\citep{liemerging}. 
The ability of retailers to offer fast and cost-effective deliveries heavily relies on making efficient dispatching and routing decisions, which involve balancing the trade-off between dispatching orders upon their arrival and waiting for future order arrivals~\citep{Klapp2018}. 
Dispatching orders early ensures that customers receive their deliveries according to their preferred delivery time, while waiting allows for delivery consolidation and scale efficiency through improved routing decisions.
However, as delivery requests are typically unknown beforehand, evaluating this trade-off requires taking into account the underlying uncertainty that arises from several factors, such as the variability of order arrival patterns, customer locations, the size of requests, and the time windows or deadlines of these requests. 
Accounting for all these sources of uncertainty in the same-day delivery environment, where decisions must be made quickly, poses considerable methodological challenges.

This paper considers the \textit{dynamic dispatch waves problem} (DDWP).
The \dynamicproblem~is a same-day delivery problem in which vehicles are dispatched at fixed decision moments (for example, every hour). 
The DDWP is found in numerous applications, such as e-commerce~\citep{Klapp2018a}, urban consolidation centers~\citep{vanHeeswijk2019}, and grocery delivery services~\citep{Liu_Luo_2022}. 
A set of new delivery requests arrive at each decision moment, and the system operator has to decide which requests to dispatch, and how to route them, while the remaining requests are postponed to future decision moments.
We consider a variant of the DDWP that was introduced in the \textit{EURO meets NeurIPS 2022} vehicle routing competition~\citep{Kool2022a}. 
In this variant, all delivery requests must be served within their requested time window and a sufficient number of vehicles is available to serve all requests on time. 
The goal is to minimize the total distance traveled.

We propose the \textit{iterative conditional dispatch} (ICD) algorithm to solve the DDWP quickly, within just a few minutes.
The algorithm leverages ideas from progressive hedging \citep{rockafellar1991scenarios} and sample scenario planning \citep{Bent2004}.
The ICD algorithm constructs a dispatching solution through the iterative resolution of sample scenarios to classify requests as dispatched, postponed or undecided.
The set of undecided requests shrinks in each iteration until a final dispatching decision is made in the last iteration.

ICD can be seen as a generalization of the dynamic stochastic hedging heuristic (DSHH) introduced by \cite{Hvattum2006} for a dynamic pickup and delivery problem with fixed decision moments.
DSHH also uses similar concepts from progressive hedging to construct a solution iteratively, and the authors use a single threshold-based mechanism that classifies which requests to serve next based on a sample scenario approach.
These decisions are then enforced by restricting the customer time windows to serve the selected requests before the next decision moment. While this adjustment is suitable for problems where vehicles are already en route and decisions concern selecting the next customers to visit, it falls short for problems involving dispatching decisions because it severely affects the quality of routes.

To overcome this issue, we use \textit{dispatch windows} to enforce dispatching decisions, leaving the customer time window and overall problem structure intact.
Moreover, we extend the single threshold mechanism of DSHH to a double threshold function: besides classifying which requests to dispatch, we also classify which requests to postpone.
Rigorous experiments on a large test bed of diverse instances demonstrate that the double threshold method applied in ICD consistently outperforms the single threshold variant by as much as 2\%.
To avoid tuning of the threshold parameter values, we also develop a non-parametric ICD variant that is as competitive as the single threshold variant.
Finally, we show that ICD performs competitively on instances from the \textit{EURO meets NeurIPS 2022} vehicle routing competition \citep{Kool2022a}, nearly matching the performance of the winning machine learning-based approach by \cite{Baty2024} that ranked first in the competition.
Although we develop and benchmark ICD for the special case of a DDWP in which all requests must be served, the idea of dispatch windows can easily be extended to more general variants of DDWPs with different objectives and constraints, such as maximizing the number of requests to be served or crowd-sourced delivery.

The outline of this paper is as follows.
Section~\ref{sec:literature} surveys the related literature. 
Section~\ref{sec:problem-description} formalizes the problem description. 
Section~\ref{sec:methods} presents the iterative conditional dispatch algorithm. 
Section~\ref{sec:experiments} describes the numerical experiments and discusses the results.
Section~\ref{sec:conclusion} concludes the paper and discusses future research ideas.

\section{Literature review}
\label{sec:literature}

The DDWP can be seen as a particular variant of the same-day delivery problem~\citep{Voccia2019, Zhang2023}.
The same-day delivery problem is a stochastic dynamic vehicle routing problem (VRP) where not all customer requests are known beforehand, but at least some arrive dynamically over the planning horizon.
The goal is to make dispatching decisions either when each customer request arrives, or at a number of fixed moments in time (for example, every hour).
Since e-commerce companies are increasingly offering same-day delivery services to enhance customer satisfaction~\citep{Klapp2020}, the rapid growth of e-commerce has resulted in a burst of recent work on same-day delivery problems:~\cite{Zhang2023} find that 90\% of same-day delivery problem-related studies have been published within the last five years. 

The most common way of modeling a stochastic and dynamic VRP is by formulating it as a Markov decision process (MDP). 
Since the state and action spaces of MDP formulations for such problems are typically very large, most studies resort to approximate algorithms. We refer to the survey of~\cite{Soeffker2022} for a classification of solution methods, and discuss here only two main ways of solving same-day delivery problems.

The first way is based on scenario sampling from a (known or approximated) stochastic process that models the dynamics of customer request arrivals.
These sampled scenarios reduce to static VRP instances, which are considerably easier to solve.
Solving each scenario in isolation, however, results in as many solutions as there are sampled scenarios.
The main difficulty in scenario sampling is how to combine the information from each scenario solution into a single, actionable decision for the original problem.

\cite{Bent2004} propose the use of a so-called \emph{consensus function} to combine the scenario solutions into a single action.
Their consensus function is based on a least-commitment approach~\citep{Stefik_1981}.
The results show a clear benefit over simply implementing the `best' scenario solution, suggesting a substantial benefit from finding consensus among the different scenarios.
\cite{Hvattum2006}~consider a dynamic pickup-and-delivery problem with fixed decision moments, and propose the DSHH to solve this problem.
The DSHH algorithm iteratively solves a set of sample scenarios, and computes the percentage of scenarios in which each request was served before the next decision moment.
They apply a threshold decision rule to these percentages to determine which request to serve.
\cite{Ghiani_Manni_Thomas_2012} consider a same-day delivery problem with a single vehicle, and propose a consensus function that selects a solution with minimal pair-wise Hamming distance from the set of scenario solutions.
\cite{Voccia2019} consider a same-day delivery problem where dispatch decisions are made as each request arrives.
Their consensus function identifies the solution that is most often observed in the scenario solutions.
\cite{Azi2012} study a same-day delivery problem with request acceptance and use a sample scenario approach to determine the relative benefit of accepting a request.
Similarly, \cite{Dayarian2020} use a sample scenario approach to determine the relative benefit of dispatching a request at the current dispatching epoch, or postponing the request for later (see \cite{Hvattum2007} for a similar idea).

Besides scenario-based approaches, there are many other ways to solve same-day delivery problems by approximating some aspect of the MDPs.
\cite{Klapp2018} study a single vehicle variant of the DDWP, in which requests are located on a line.
They formulate a rollout policy and an approximate linear programming-based dynamic policy to solve the problem.
\cite{Klapp2018a} generalize this setting to general network topologies.
\cite{vanHeeswijk2019} formulate an MDP for a DDWP, and solve the MDP using approximate dynamic programming with a value function approximation through basis functions.
They test their method on instances with up to forty customer requests arriving between decision moments, and find that they outperform several myopic policies by 6\% to 20\%.
\cite{Ulmer2019} formulate an MDP for a same-day delivery problem in which vehicles are allowed to return to the depot before all assigned requests have been served.
They propose a value function approximation and show that preemptive depot returns can substantially increase the number of deliveries.
\cite{Ulmer_2020} integrates pricing and routing for same-day delivery in an e-commerce setting.
They use a value function approximation based on state space aggregation to solve the problem efficiently.
\cite{Liu_Luo_2022} present a dynamic programming formulation and structured approximation method for a DDWP.
They prove several structural results on the quality of the resulting solution, and validate the approach numerically through a real-world case in grocery and food delivery. 
Recently~\cite{Heitmann2023} combined the two approaches in a dial-a-ride setting.
They use value function approximation to efficiently decide whether to accept a request, while using scenario sampling and the consensus function of~\cite{Bent2004} to decide on the routing problem.

Scenario sampling and other MDP-based techniques both have to work around the issue that the underlying vehicle routing problem in same-day delivery problems is, in general, difficult to solve.
This makes accurately estimating the (expected) cost of a dispatching decision a difficult task, since that typically requires solving many routing problems in a decision tree.
This is not an option in the time-constrained setting of the~\dynamicproblem~we consider here.
Several techniques have been developed to curtail the computational effort needed to explore the solution space.
Dominance checks can be used to eliminate sub-optimal nodes in the decision tree without sacrificing optimality~\citep{goodson2016restocking}.
Alternatively,~\cite{secomandi2009reoptimization} restrict the state space heuristically and~\cite{zhang2022offline} impose a-priori rules that restrict the customer search space based on the current solution.
In particular, they only allow vehicles to visit subsequent customers within a predefined radius. 
\cite{Hvattum2006} use principles from progressive hedging~\citep{rockafellar1991scenarios} to iteratively fix decisions for some requests.
This progressively builds a solution of high quality, which can be much faster than exhaustively solving the scenario sampling or MDP problem.

\section{Problem description}
\label{sec:problem-description}

In this section, we formulate the \dynamicproblem~variant that was presented in the \textit{EURO meets NeurIPS 2022} vehicle routing competition \citep{Kool2022a}.
We introduce notation and provide a general problem outline in Section~\ref{subsec:outline}.
Thereafter, we introduce the vehicle routing problem with time windows and dispatch windows (VRPTW-DW) in Section~\ref{subsec:static_vrptwdw}.
Finally, in Section~\ref{sec:mdp}, we formulate the \dynamicproblem~as an MDP.

\subsection{General outline}
\label{subsec:outline}
We consider a planning horizon of $H$ time units, which is further divided into a set of decision epochs $\mathcal{T} =$ $\{ 1, 2, \dots, t_{\text{final}} \}$, where $t_\text{final}$ represents the final epoch.
We assume each epoch $t \in \mathcal{T}$ starts at time $T_t \ge 0$, with $T_t > T_{t - 1}$ for $t > 1$.

A set of new delivery requests $\omega_t$ is revealed at the start of each decision epoch $t \in \mathcal{T}$.
The set of requests $\omega_t$ has support $\Omega_{t}$.
We assume the support $\Omega_t$ is known, either through direct knowledge or by approximation from historical data. 
Each request $i \in \omega_t$ requires a service time of $\mu_i \ge 0$ time units when serviced at the request location. 
Service may only start between the request's time window $\left[e_i, l_i \right]$, with $T_t \le e_i < l_i \le H$. 
Early arrival at a request location is allowed, but service must wait until the time window opens.
The demand of a request is denoted by $q_i \ge 0$ units of load capacity of vehicles.
The release time of a request is denoted by $r_i = T_t$.
The travel time between the locations of two requests $i$ and $j$ is denoted by $\tau_{ij} \ge 0$, and the cost of travel is denoted by $c_{ij} \ge 0$.
Deliveries are made with a fleet of homogeneous vehicles $\mathcal{K}$, where each vehicle has a capacity of $Q > 0$ units.
A feasible route dispatched in epoch $t$ starts at the depot at time $T_t$ and must return to the depot before the end of the planning horizon $H$.
The sum of the requested demands on the route may not exceed the vehicle capacity.

We consider that operators always have access to sufficient vehicles at the start of each decision epoch.
This assumption is different from the common restriction in same-day delivery problems that only a limited number of vehicles is available, while allowing for some requests not to be served~\citep{Voccia2019, Klapp2018}. Despite these differences, our problem and solution approach capture the dynamics between waiting and dispatching decisions to improve operational efficiency, which is one of the core challenges of same-day delivery problems. To offer a more realistic view, we include an alternative model in Appendix~\ref{app:heterogeneous-fleet-extension} with limited fleet availability and accompanying numerical results.

At the start of each decision epoch, the operator must decide which of the known requests to dispatch and which ones to postpone to consolidate with new requests that may arrive in later epochs.
The objective is to dispatch requests and route vehicles such that all requests are served within their time windows, and to minimize the total travel cost of the routes in current and future epochs.
This is formalized in Section~\ref{sec:mdp}. 

\subsection{Static vehicle routing problem with dispatch windows}
\label{subsec:static_vrptwdw}
Before the dynamic problem is formulated using MDP terminology, we present the VRPTW with \textit{dispatch windows} (VRPTW-DW) as a preliminary problem.
The VRPTW-DW is the static and deterministic version of the DDWP. Similar models have been presented before in \cite{Klapp2018a}~and~\cite{vanHeeswijk2019}.
Let $\mathcal{N}$ be the index set of delivery requests, and let $T \ge 0$ denote the earliest time that a request may be dispatched.
The depot is represented by two indices $\sdepot$ and $\edepot$, which correspond to the start and end of a route, respectively.
The depot has no service time ($\mu_{\sdepot} = \mu_{\edepot} = 0$), no demand ($q_{\sdepot}=q_{\edepot}=0$), and a time window that spans the planning horizon ($\left[e_{\sdepot}, l_{\sdepot} \right] =\left[e_{\edepot}, l_{\edepot} \right]= \left[0, H \right]$).
Let $\mathcal{V} = \{\sdepot, \edepot \} \cup \mathcal{N}$ define the set of locations.
A dispatch window $[r^-_i, r^+_i]$ is a time interval during which a route containing request $i \in \mathcal{N}$ must leave the depot.
Unless otherwise specified, we assume $r^-_i = r_i$ and $r^+_i = H$ for each request $i \in \mathcal{N}$.
Two sets of decision variables are used to define a solution.
First, $x_{ijk} \in \{0, 1\}$ tracks the routing decisions: it is $1$ if vehicle $k \in \mathcal{K}$ travels from location $i \in \mathcal{V}$ to $j \in \mathcal{V}$, and $0$ otherwise.
Second, $\theta_{ik} \ge 0$ determines the time at which service starts at location $i \in \mathcal{V}$ by vehicle $k \in \mathcal{K}$.

Let
\[ c(x, \theta) = \sum_{k \in \mathcal{K}} \sum_{i, j \in \mathcal{V}} c_{ij} x_{ijk} \]
denote the total (travel) cost of a decision pair $(x, \theta)$.
A formulation for the VRPTW-DW is given by:
\begin{subequations}
\label{milp:vrptw-dw}
\begin{align}
   \text{VRPTW-DW}(\mathcal{N}, T): \min_{x,\theta} \quad & c(x, \theta) \label{con:objective}
    \\ \text{s.t. } \quad 
    & \sum_{k \in \mathcal{K}} \sum_{j \in \mathcal{V}} x_{i j k}=1 & \forall i \in \mathcal{N}, \label{con:every-customer-assigned}
    \\ & \sum_{j \in  \mathcal{V}} x_{\sdepot j k}=1 & \forall k \in  \mathcal{K}, \label{con:start-at-depot}
    \\ & \sum_{i \in  \mathcal{V}} x_{i, \edepot, k}=1 & \forall k \in \mathcal{K}, \label{con:end-at-depot}
    \\ & \sum_{j \in \mathcal{N}} x_{i j k} = \sum_{j \in  \mathcal{N}} x_{j i k} & \forall k \in \mathcal{K}, i \in  \mathcal{N}, \label{con:flow-conservation}
    \\ & x_{i j k}(\theta_{ik}+\mu_{i}+\tau_{ij} - \theta_{jk}) \leq 0  & \forall k \in \mathcal{K}, i, j \in  \mathcal{V}, \label{con:arrival-times} 
    \\ & e_{i} \leq \theta_{i k} \leq l_{i} & \forall k \in \mathcal{K}, i \in \mathcal{V}, \label{con:time-windows}
    \\ & \max\{r^-_i, r^-_j\} x_{ijk} \le \theta_{0k}  &  \forall k \in \mathcal{K}, i, j \in  \mathcal{N}, \label{con:dispatch-window-min}
    \\ & \theta_{0k} \le \min\{r^+_i, r^+_j\} x_{ijk}  &  \forall k \in \mathcal{K}, i, j \in  \mathcal{N}, \label{con:dispatch-window-max}
    \\ & \theta_{0k} \geq T & \forall k \in \mathcal{K}  \label{con:earliest-departure}
    \\ & \sum_{i \in \mathcal{V}} q_{i} \sum_{j \in \mathcal{V}} x_{i j k} \leq Q & \forall k \in  \mathcal{K},   \label{con:capacity}
    \\ & x_{i j k} \in\{0,1\} & \forall k \in  \mathcal{K}, i, j \in  \mathcal{V},
    \\ & \theta_{ik} \geq 0 & \forall k \in  \mathcal{K}, i \in  \mathcal{V}. \label{eq:constraint-1vrptw-rt}
\end{align}
\end{subequations}
The objective~\eqref{con:objective} minimizes the total travel cost. 
Constraints~\eqref{con:every-customer-assigned} ensure that each request is assigned to exactly one vehicle. 
Constraints~\eqref{con:start-at-depot} and \eqref{con:end-at-depot} define that each vehicle starts and ends at the depot, respectively. 
Constraints~\eqref{con:flow-conservation} ensure flow conservation at each request location.
Constraints~\eqref{con:arrival-times} and \eqref{con:time-windows} guarantee feasibility with respect to the arrival times and time windows, respectively. 
Constraints~\eqref{con:dispatch-window-min} and \eqref{con:dispatch-window-max} ensure feasibility with respect to the request dispatch windows, and constraints~\eqref{con:earliest-departure} ensure that no vehicle is dispatched before time $T$.
Finally, constraints~\eqref{con:capacity} ensure feasibility with respect to vehicle capacity.

The problem formally described by~\eqref{milp:vrptw-dw} is relevant to the dynamic problem in three ways.
First, and most importantly, it is repeatedly solved heuristically by our solution approach in Section~\ref{sec:methods} to determine dispatching actions in each epoch.
Second, it is used to evaluate the cost of a dispatching decision at each decision moment (see Section~\ref{sec:contribution-and-objective}). 
Third, with perfect information, the \dynamicproblem~reduces to a static and deterministic VRPTW-DW problem.
Let $\mathcal{N}_{t}= \bigcup_{t' = 1}^t \omega_{t'}$ denote the set of all realized requests up to epoch $t$.
Solving VRPTW-DW($\mathcal{N}_{t_{\text{final}}}, T_1)$ now provides a lower bound on the cost of the dynamic counterpart for this realization.
We refer to this large static problem as the \textit{hindsight} problem and use it in Section~\ref{sec:experiments} for benchmarking.

\subsection{Markov decision process formulation}
\label{sec:mdp}
We now formulate the \dynamicproblem~as an MDP.

\subsubsection{State space}
The state space specifies the current state of the system in reference to requests that are known but not yet dispatched.
The state space at epoch $t$ can be simply represented as a set $s_t = \{i\in \mathcal{N}_{t}~\mid~\text{request $i$ has not been dispatched}\}$.

\subsubsection{Action space} \label{S:actionSpace}
Given the current state $s_t$ in epoch $t$, an action $a_t$ is a subset of available requests that are chosen to be dispatched in epoch $t$, that is, $a_t \subseteq s_t$.
Let $m_t \subseteq s_t$ denote the set of requests in epoch $t$ that must be dispatched because of feasibility requirements. 
Specifically, we define
\[ m_t = \left\{ i \in s_t~\mid~\text{VRPTW-DW($\{i\}$, $T_{t+1}$) is infeasible} \right\} \]
as the set of requests that cannot be feasibly dispatched as a single request in the next epoch at time $T_{t+1}$. 
The complete action space $\mathcal{A}(s_t)$ is then specified as
\[ \mathcal{A}(s_t) = \left\{ a_t \subseteq s_t ~\mid~m_t \subseteq a_t \right\}. \]
At its largest, the action space can be equal in size to $2^{\vert s_t \vert}$, the power set of the current state $s_t$.

Illustrative of the literature, \cite{Voccia2019} and \cite{vanHeeswijk2019}, among others, use more involved state space representations where route information such as return times are also specified.
Also, \cite{Ulmer2020} argue for the use of a route-based representation of state space.
We chose to simplify our representation for two reasons.
First, unlike the aforementioned studies, our assumption of unlimited available vehicles ensures that the dispatching of vehicles in the current epoch does not impact decisions in later epochs.
Therefore, including additional route-specific information would only serve to expand the state space.
Second, the existence of multiple optimal routes with the same cost would increase the size of the action space.
Our set representation avoids this symmetry issue and substantially reduces the size of the action space.

\subsubsection{State transitions} 
The transition to the next state $s_{t+1}$ depends on the chosen action $a_t$ and the unknown future realization $\omega_{t+1}$.
For each realization $\omega_{t+1}$ from the support $\Omega_{t+1}$, the next epoch's state is given by $s_{t+1}= (s_t \setminus a_t) \cup \omega_{t+1}$. 

\subsubsection{Optimality conditions}
\label{sec:contribution-and-objective}
The direct cost $C(s_t, a_t)$ of a specific action $a_t$ in state $s_t$ is given by the cost of routing the requests in $a_t$ at time $T_t$.
Let $(x^*, \theta^*)$ be the optimal solution to VRPTW-DW($a_t, T_t$).
The direct cost is then given by $C(s_t, a_t) = c(x^*, \theta^*)$.
The objective of the DDWP is to select for each epoch $t \in \mathcal{T}$ a minimum cost action, such that the following (Bellman) optimality conditions are satisfied:
\begin{align}
    \label{eq:optimalPolicyEquation}
    V(s_t) = 
    \begin{dcases}
        C(s_{t}, s_{t}) & \text{if } t = |\mathcal{T}|, \\
        \min_{a \in \mathcal{A}(s_t)} \big[ C(s_t, a)  + \mathbb{E}_{\omega_{t + 1}} [ V(s_{t+1}) ] \big] &  \text{otherwise}.
    \end{dcases}
\end{align}
The conditions of~\eqref{eq:optimalPolicyEquation} lend themselves naturally to a dynamic programming-based solution approach.
However, there are three main difficulties in pursuing a solution approach in this direction.
First, the number of actions $\mathcal{A}(s_t)$ is exponential in the size of the state $s_t$. 
Second, obtaining the costs $C(s_t, a_t)$ for a single state and action pair $(s_t, a_t)$ requires solving the challenging VRPTW-DW of Section~\ref{subsec:static_vrptwdw}.
Third, an accurate specification of the transition probabilities that are needed to evaluate the expected future costs are often not available in practice.
In the next section, we propose an algorithm that is tailored to overcome these challenges.

\section{The iterative conditional dispatch algorithm}
\label{sec:methods}

This section presents the \textit{iterative conditional dispatch} (ICD) algorithm to solve the \dynamicproblem.
The ICD algorithm is an iterative solution construction procedure based on a sample scenario approach to tractably approximate the exponential state and action spaces outlined in Section~\ref{sec:mdp}.
The pseudo-code of the ICD algorithm is given in Algorithm~\ref{alg:icd} for a fixed epoch $t \in \mathcal{T}$.
The algorithm maintains two sets of decisions: the set of dispatched requests $d_t$, and the set of postponed requests $p_t$.
At initialization, we initially set $d_t = m_t$ and $p_t = \emptyset$.
Then, the ICD algorithm iteratively applies two steps: step I and step II.
In step I (Section~\ref{subsec:solve_sample_scenarios}), we first generate a set of sample scenarios.
These sample scenarios each represent one potential realization of future requests.
Then, with the previously fixed decisions given by $d_t$ and $p_t$, we solve the sample scenarios as static VRPTW-DWs (Section~\ref{subsec:static_vrptwdw}) to identify the dispatch and postponement decisions that are made for each sample instance. 
In step II (Section~\ref{subsec:consensus_functions}), we use the scenario solutions to classify requests that are not in $d_t \cup p_t$ into either the dispatch set $d_t$, the postponement set $p_t$, or leave the request undecided.
Note that we maintain the invariant $d_t \cap p_t = \emptyset$, so a request cannot be classified as both dispatched and postponed.
The classification of step II reduces the number of undecided requests in subsequent iterations.
Upon finishing a specified number of iterations or when all requests have been classified (that is, when $s_t = d_t \cup p_t$), the dispatch action $a_t$ is given by $d_t$.
Sections~\ref{subsec:solve_sample_scenarios} and~\ref{subsec:consensus_functions} describe steps I and II in more detail.
\begin{algorithm}
\caption{Iterative conditional dispatch in epoch $t \in \mathcal{T}$.}
\label{alg:icd}
\begin{algorithmic}[1]
    \State Initialize $d_t = m_t$ and $p_t = \emptyset$
    \While{number of iterations not exceeded \textbf{and} $s_t \neq d_t \cup p_t$}                    
        \State \textbf{Step I}: generate scenarios $\mathcal{S}$ and solve them as VRPTW-DWs conditioned on $d_t$ and $p_t$
        \State \textbf{Step II}: update $d_t$ and $p_t$ based on scenario solutions
    \EndWhile
    \State $a_t = d_t$
\end{algorithmic}
\end{algorithm}

\subsection{Step I: solving sample scenarios}
\label{subsec:solve_sample_scenarios}

Step I requires the generation of $|\mathcal{S}|$ sample scenarios and solving each of the resulting static routing problems. 
Formally, a single scenario $S$ is constructed as follows. 
A scenario $S \in \mathcal{S}$ is a set of requests formed by the current epoch state $s_t$ and a sample path realization $(\tilde{\omega}_{t+1}, \tilde{\omega}_{t+2}, \dots, \tilde{\omega}_{t+L})$, where $\tilde{\omega}_{t'}$ sampled from $\Omega_{t'}$ denotes the set of future requests that arrive in decision epoch $t'$ and $L$ denotes the number of lookahead epochs.
Thus,
\[ S = s_t \cup \bigcup_{t' = t + 1}^{t + L} \tilde{\omega}_{t'}. \]

Each scenario is solved as the problem defined by $\text{VRPTW-DW}(S, T_t)$.
To incorporate the dispatch and postponement decisions $d_t$ and $p_t$, we make the following adjustments to the dispatch windows $\left[ r^-_i, r^+_i \right]$ for each request $i \in S$.
For dispatched requests $i \in d_t$, the dispatch window is set to $\left[r^-_i, r^+_i\right] = \left[T_t, T_t \right]$, forcing these requests to be dispatched in the current epoch $t$ at decision moment $T_t$.
For postponed requests $i \in p_t$, the dispatch window is set to $\left[r^-_i, r^+_i\right] = \left[T_{t+1}, H \right]$, where $T_{t+1}$ is the decision moment of the next epoch $t + 1$.
All other requests $i \in S \setminus (d_t \cup p_t)$ retain their default dispatch windows.

A feasible solution $(x, \theta)$ to VRPTW-DW($S, T_t)$ is readily mapped back to a \textit{scenario solution} $a_S$.
In particular, we can infer the request and vehicle pairings from the $x$ decisions.
Then, the epoch in which each request was dispatched can be read from the route starts $\theta_{0k}$ for each vehicle $k \in \mathcal{K}$.
Formally, a request $i \in s_t$ is dispatched in epoch $t$ if $\theta_{0k_i} = T_t$, where $k_i$ is the vehicle to which $i$ is assigned.
We thus have
\[ a_S = \{ i \in s_t~\mid~\theta_{0k_i} = T_t \}. \]

One of the challenges in sample scenario approaches is that the scenario problems are often hard to solve optimally.
We show in the numerical experiments (Section~\ref{sec:experiments}) that it is not necessary to solve these problems optimally to obtain good solutions using ICD.

\subsection{Step II: updating \texorpdfstring{$d_t$}{dispatch} and \texorpdfstring{$p_t$}{postponement} using consensus functions}
\label{subsec:consensus_functions}
Step II involves classifying requests based on the scenario solutions to update the action sets $d_t$ and $p_t$.
We use consensus functions to make the classification.
Our first consensus function is based on thresholds, and our second on the Hamming distance.

\subsubsection{Threshold-based consensus.}\label{sec:threshold}
The threshold-based consensus function classifies each request based on how often it is dispatched in the scenario solutions.
First, we define the dispatch score $0 \le \phi(i) \le 1$ of a known request $i \in s_t$ as the relative frequency of it being dispatched among the sample scenario solutions.
Formally, the dispatch score $\phi(i)$ of a request $i \in s_t$ is given by
\begin{align*}
    \phi(i) & = \frac{1}{\vert \mathcal{S} \vert} \sum_{S \in \mathcal{S}} \mathbbm{1}_{i \in a_S},
\end{align*}
where $\mathbbm{1}_{i \in a_S}$ is an indicator that is $1$ if $i$ is dispatched by scenario solution $a_S$, and $0$ otherwise.

The threshold-based consensus function is characterized by two separate threshold parameters $\epsilon_D$ and $\epsilon_P$ that determine the cut-off points for the assignment of requests to action sets $d_t$ and $p_t$.
Given the dispatch scores $\phi(i)$, we update $d_t$ and $p_t$ as follows:
\begin{align*}
    d_t &= d_t \cup \{n \in s_t~\mid~ \phi(i) \geq \epsilon_D \}, \\ 
    p_t &= p_t \cup \{n \in s_t~\mid~ \phi(i) < \epsilon_P \}.
\end{align*}

The dispatch score $\phi(i)$ can be interpreted as a measure of the probability that a request $i$ should be dispatched in epoch $t$. 
If the dispatch score of a request is high, then the request is a good candidate for dispatching. 
Conversely, requests with high postponement probabilities ($1-\phi(i)$) should be postponed. 
The threshold parameters $\epsilon_D$ and $\epsilon_P$ reflect the tolerance level for classifying a request as dispatched or postponed, respectively.
Requests with a dispatch score in $\left[\epsilon_P, \epsilon_D\right)$ remain unclassified in the current iteration, as their dispatch probabilities are not high enough to be classified as dispatched, and not low enough to be classified as postponed.
In later iterations, the dispatch scores of the unclassified requests are expected to change due to the updated action sets $d_t$ and $p_t$.

Note that the disjoint property of $d_t$ and $p_t$ requires $\epsilon_D \ge \epsilon_P$, since otherwise a request can be marked for dispatch and postponement simultaneously, which violates the feasibility of the action.
Also observe that when $\epsilon_D + \epsilon_P = 1$, the ICD reduces to a non-iterative procedure as all requests are marked in one iteration.

\subsubsection{Similarity-based consensus.}\label{sec:hamming}
Determining the optimal threshold values for the threshold-based consensus function requires substantial parameter tuning. 
Therefore, we develop an additional parameter-free consensus function that selects a scenario solution with the lowest average Hamming distance to the other scenario solutions.
The Hamming distance counts the number of requests that two solutions decide on differently in their dispatching actions.
We update the partial dispatch set by selecting the scenario $S^*$ whose solution $a_{S^*}$ obtains the minimum average Hamming distance:
\begin{align*}
      S^* &= \arg \min_{S \in \mathcal{S}} \frac{1}{|\mathcal{S}|} \sum_{S' \in \mathcal{S}} | a_S \setminus a_{S'} | + | a_{S'} 
      \setminus a_S | , \\
      d_t &= d_t \cup a_{S^*}.
\end{align*}
The intuition behind this particular update is as follows.
When the solution that is most similar to all other scenario solutions dispatches a request, it is likely that most other solutions also made this decision, and dispatching is a good idea.
This solution may also include requests that are not as frequently dispatched, as this may have been a good choice to balance out the direct cost and the future costs for the specific scenario.
Thus, when fixing a postponement decision we are more conservative: only when \textit{all} scenario solutions postpone the request then we classify it as postponed.
The partial postponement set is hence updated as:
\[ p_t = p_t \cup \{i \in s_t~\vert~i \not\in a_S \text{ for all } S \in \mathcal{S} \}. \]

\section{Numerical experiments}\label{sec:experiments}
This section presents the computational experiments that show the effectiveness of our ICD algorithm.
All numerical experiments are conducted on a single thread of an AMD EPYC 7H12 CPU.
Our code is implemented in Python and is released at \url{https://github.com/leonlan/dynamic-dispatch-waves} under the MIT license.
The static VRPTW-DW problems are solved using PyVRP version 0.6.3, a high-performance and open-source heuristic VRP solver implemented in C++ and Python~\citep{wouda_pyvrp_2023}. It builds upon the hybrid genetic search algorithm of~\cite{vidal2013hgs} and~\cite{vidal2022hgscvrp}.
Specifically for this paper, we extend PyVRP with support for dispatch windows to solve the VRPTW-DW of Section~\ref{subsec:static_vrptwdw}.
Implementation details of PyVRP and its dispatch windows extension are given in Appendix~\ref{app:implementation-details}.

The rest of this section is structured as follows.
We describe test instances in Section~\ref{sec:instances}.
Section~\ref{sec:benchmark-algorithms} and Section~\ref{sec:parameter-settings} describe the benchmark algorithms and their parameter settings, respectively.
Section~\ref{sec:results} presents the numerical results for our test instances and investigates the individual components of ICD.
Finally, in Section~\ref{sec:euro-neurips}, we evaluate ICD on instances from the \textit{EURO meets NeurIPS 2022} vehicle routing competition.

\subsection{Benchmark instances}\label{sec:instances}
Benchmark instances for the DDWP are based on the topologies of the extended Solomon instances~\citep{Homberger1999}.
We replace the existing time windows of these instances with a new time window generation process to better control the stochastic dynamics of the benchmark instances. 
We consider all of the three location dispersion sets of the Solomon instances: randomly distributed (\textit{R}), clustered (\textit{C}), and randomly distributed and clustered (\textit{RC}). 
Two large static instances with 1\,000 customers from each of these three categories are considered.
Specifically, we select \textit{R1\_10\_1, R2\_10\_1} for \textit{R}, \textit{C1\_10\_1, C2\_10\_1} for \textit{C}, and \textit{RC1\_10\_1, RC2\_10\_1} for \textit{RC}.
The Euclidean distances between each combination of two locations in these static instances are normalized such that the furthest location can be served within at most one epoch duration. 
Moreover, we assume that the travel cost and duration between locations are equal to the travel distance.

\begin{figure}
  \centering
  \includegraphics[width=\linewidth]{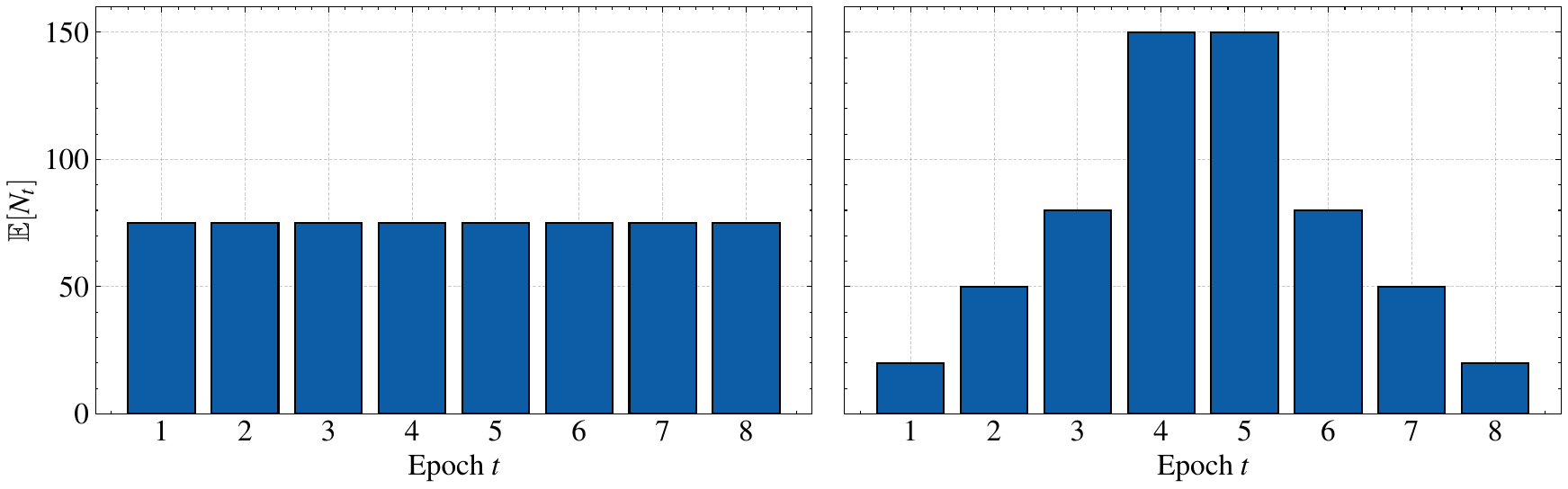}
  \caption{
    The expected number of arrivals per epoch for the homogeneous (left) and unimodal (right) arrival processes.
  }
  \label{fig:arrivals}
\end{figure}

We now explain how the dynamic instances for the DDWP are created. 
The planning horizon of $H = 8$ hours is discretized into eight epochs $\mathcal{T} = \{1, \dots, 8\}$ of equal duration (1 hour each).
At each epoch $t$, we sample $N_t$ requests and their attributes.
To adequately capture the dynamics of different arrival processes, two types of arrival processes $\{N_t\}_{t \in \mathcal{T}}$ are considered.
We call these two arrival processes \textit{homogeneous} and \textit{unimodal}, respectively, and illustrate them graphically in Figure~\ref{fig:arrivals}.
The number of arrivals in both arrival processes follows a discrete uniform distribution with bounded support on $\left[ \lfloor 0.9 \mathbb{E}[N_t] \rfloor, \lfloor 1.1 \mathbb{E}[N_t] \rfloor \right]$.
For the homogeneous arrival process, the expected number of arrivals remains constant at $\mathbb{E}[N_t] = 75$ requests for every epoch $t$ in $\mathcal{T}$. 
The unimodal process has epoch-dependent arrival rates with a peak in the middle of the planning horizon. 
The expected arrivals in the eight epochs are $(20, 50, 80, 150, 150, 80, 50, 20)$ delivery requests, where the $t^{\text{th}}$ element in the sequence denotes $\mathbb{E}[N_t]$.

\begin{figure}
  \centering
  \includegraphics[width=\linewidth]{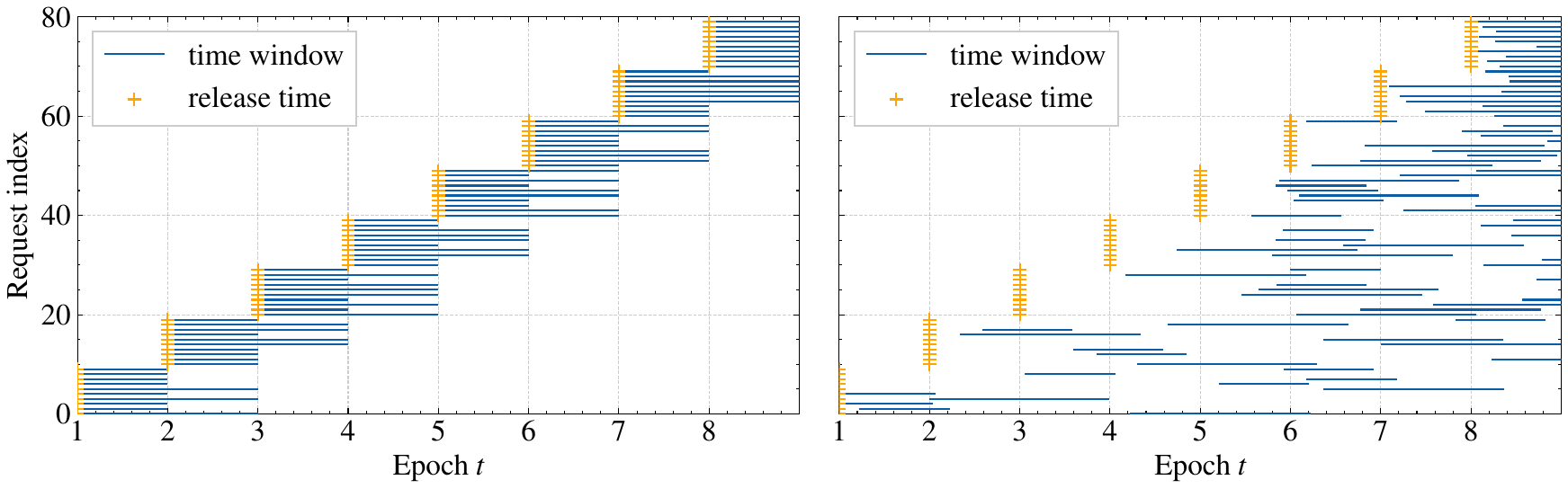}
  \caption{
    Two dynamic instances with ten requests per epoch and maximum time window width of two.
    The left figure shows an instance with deadlines and the right figure shows an instance with regular time windows.
  }
  \label{fig:tws}
\end{figure}

Request attributes corresponding to location, demand and service time are sampled uniformly and independently from the static instance. 
For the time windows, we examine six variants, which are characterized by \textit{width} and \textit{type}.
We consider three different maximum time window widths: two hours, four hours, and eight hours.
For every request, the width of the time window is sampled from a discrete uniform distribution between one and the maximum time window width.
Let $W_i$ denote the width of the time window for request $i$.
We examine two types of time windows: deadlines and regular time windows. 
Deadlines are defined such that the time window starts at the request's release time $r_i$ and ends at $r_i + W_i$.
Regular time windows begin at a time point uniformly sampled between the request's release time $r_i$ and the planning horizon $H$. 
They end at this starting time plus $W_i$.
Finally, all time windows are truncated to end no later than the planning horizon.
We label the resulting variants as \textit{DL2, DL4, DL8, TW2, TW4}, and \textit{TW8}, where \textit{DL} represents deadlines, \textit{TW} represents regular time windows, and the digit indicates the maximum time window width in hours. 
The time window variants DL2 and TW2 are visualized in Figure~\ref{fig:tws}. 

In total, we have three topologies, two arrival processes, and six time window variants.
We use a full factorial design of these different parameters, and consider all 36 combinations of these parameters. 

\subsection{Benchmark algorithms}\label{sec:benchmark-algorithms}
We implement several variants of the ICD algorithm.
The first variant uses the double threshold consensus function of Section~\ref{sec:threshold}, which we refer to as \textit{ICD-double}. 
The second variant uses the Hamming distance consensus function of Section~\ref{sec:hamming} and we refer to this variant as \textit{ICD-Hamming}.
We also consider two variants that use only a single threshold, which allows us to study the marginal contribution of both dispatch and postponement thresholds.
The ICD variant with only a dispatch threshold closely resembles the DSHH of~\cite{Hvattum2006}, and we similarly refer to this variant as \textit{DSHH}.
The other single threshold variant uses only a postponement threshold, and we refer to this variant as \textit{ICD-postpone}.
Because the ICD-postpone variant does not explicitly classify dispatch requests, the final action $a_t$ is determined by the requests that are not postponed, that is, $a_t = s_t \setminus p_t$.
Finally, we also implement a rolling horizon (RH) policy for comparison. The RH policy only samples a single scenario and uses the solution to this instance to determine which requests to dispatch. Note that the RH policy can be interpreted as a variant of ICD-Hamming using just a single sample scenario and one iteration.
All algorithms are compared against the cost of the solution obtained from heuristically solving the hindsight problem using PyVRP (see Section~\ref{subsec:static_vrptwdw}).

\subsection{Parameter settings}\label{sec:parameter-settings}
All of the ICD variants employ three iterations, thirty scenarios per iteration, and one lookahead epoch.
The ICD parameters were selected based on manual testing and observations from the related literature~\citep{Hvattum2006,Voccia2019}.
To find the best threshold parameters, we evaluated the threshold-based ICD variants on 36 dynamic instances featuring homogeneous arrival processes and we generated five different request realizations for each. 
Four separate values for both the dispatch threshold (0.4, 0.5, 0.6, 0.7) and the postponement threshold (0.2, 0.3, 0.4, 0.5) were tested. 
ICD-double was evaluated on a full factorial design of these threshold values, while DSHH and ICD-postpone were evaluated using only the dispatch and postponement values, respectively.
The best results for ICD-double were achieved with dispatch and postponement thresholds of $\epsilon_D = 0.5$ and $\epsilon_P = 0.2$, respectively. 
DSHH obtained the best results with a dispatch threshold of $\epsilon_D = 0.5$, and ICD-postpone obtained the best results with a postponement threshold of $\epsilon_P = 0.3$.
For the RH policy, we use one lookahead epoch.

The total solving time in a single epoch is set to just 120 seconds.
We keep this time limited to ensure our methods can be applied in practical same-day delivery settings, where dispatching decisions also need to be made quickly.
In each epoch, the ICD methods solve 90 scenarios in total, to which we assign a time budget of 90 seconds, meaning that ICD solves each scenario using a time limit of only one second. 
Then, an additional 30 seconds is spent on computing a route plan and associated cost for the selected dispatching action.
Although our implementation does not leverage parallelization, it is important to note that the scenario-solving process could be easily parallelized to substantially reduce the total solving time.
The RH policy uses the full 90 seconds to solve one scenario and 30 seconds to compute the route plan.
The hindsight solution is obtained by solving the hindsight problem for 600 seconds.

An important element of the ICD framework is to solve the static scenarios as efficiently as possible.
We refer to Appendix~\ref{sec:hgs-parameter-tuning} for the parameter tuning procedure of the static solver to solve scenarios well.

\subsection{Results on the benchmark instances}
\label{sec:results}
In this section, we evaluate performance on the 36 instance classes described in Section~\ref{sec:instances}.
For each instance class, we generate 200 different dynamic instances (100 per selected static instance in each of the R, C, and RC sets).
This results in a total of 7\,200 instances.
Each method is run on all 7\,200 instances and we record the final cost of each such run. 
Table~\ref{tab:main-results} presents an overview of the results. 
It presents the average percentage gap to the cost of the hindsight solution, for each method and instance class.
Entries in bold indicate the best average gaps for each instance class (row). 
Furthermore, an asterisk (*) next to the best result indicates that the superior performance of the method is significant at the 0.05 level compared to all other methods in pairwise t-tests (with Bonferroni correction).
A complete table of the statistical test results can be found in Table~\ref{tab:main-results-p-values} (Appendix~\ref{app:results}).

Table~\ref{tab:main-results} shows that among all instance classes, ICD-double solutions have an average gap of 8.05\% compared to the hindsight solutions.
On average, ICD-double has a lower gap of 1.5\%, 0.7\%, 1.6\% and 0.8\% compared to RH, DSHH, ICD-postpone and ICD-Hamming, respectively.
When we consider results by instance class (an individual row in Table~\ref{tab:main-results}), ICD-Hamming and ICD-double are better on average in 33 out of 36 instance classes, and ICD-double's improvement is statistically significant for 14 out of 36 instance classes. 
Deadline instances are arguably easier and allow for less significant savings by correctly anticipating future orders, as evidenced by the lower gaps of all methods on these instances compared to the hindsight solutions.

\begin{table}[ht]
\caption{
    Overview of results.
    The values represent the average percentage gap with respect to the hindsight solution over 200 different instances with 600 expected total number of requests.
    The best average gap in each row is marked in bold.
    An asterisk (*) indicates the method is statistically significantly better than the other methods for that instance at the 0.05 level (with Bonferroni correction).}
\label{tab:main-results}
\small
\centering
\begin{tabular}{llllrrrrr}
\toprule
\multicolumn{3}{c}{Instance} && \multicolumn{5}{c}{Methods} \\
\cmidrule{1-3} \cmidrule{5-9}
Topology & Arrival & TW &&  RH &  DSHH &  ICD-postpone &  ICD-Hamming &  ICD-double \\
\midrule
\multirow[t]{12}{*}{R} & \multirow[t]{6}{*}{HOM} & DL2 &  & 4.02 & 2.66 & 3.07 & \textbf{2.46} & 2.57 \\
 &  & DL4 &  & 8.66 & 6.84 & 8.00 & 6.38 & \textbf{6.35} \\
 &  & DL8 &  & 11.74 & 10.20 & 12.35 & 9.86 & \textbf{9.73} \\
     \cmidrule(lr){3-9}
 &  & TW2 &  & 13.38 & 12.40 & 13.47 & 12.64 & \textbf{11.02*} \\
 &  & TW4 &  & 13.84 & 13.07 & 14.29 & 13.49 & \textbf{11.76*} \\
 &  & TW8 &  & 13.52 & 13.17 & 14.10 & 13.54 & \textbf{11.66*} \\
     \cmidrule(lr){2-9}
 & \multirow[t]{6}{*}{UNI} & DL2 &  & 4.74 & 2.73 & 3.01 & 2.78 & \textbf{2.69} \\
 &  & DL4 &  & 10.20 & 8.03 & 8.87 & 7.79 & \textbf{7.56} \\
 &  & DL8 &  & 13.55 & 12.13 & 13.17 & 12.19 & \textbf{11.41*} \\
     \cmidrule(lr){3-9}
 &  & TW2 &  & 14.15 & 15.12 & 14.91 & 15.39 & \textbf{13.40*} \\
 &  & TW4 &  & 14.77 & 15.86 & 15.91 & 16.48 & \textbf{14.46} \\
 &  & TW8 &  & 14.63 & 16.24 & 16.20 & 16.90 & \textbf{14.56} \\
 \midrule
\multirow[t]{12}{*}{C} & \multirow[t]{6}{*}{HOM} & DL2 &  & 1.80 & \textbf{0.12} & 0.16 & 0.15 & 0.12 \\
 &  & DL4 &  & 3.00 & 1.41 & 1.91 & \textbf{1.29} & 1.38 \\
 &  & DL8 &  & 4.28 & 3.12 & 4.45 & \textbf{2.88} & 3.08 \\
     \cmidrule(lr){3-9}
 &  & TW2 &  & 7.68 & 6.39 & 7.99 & 6.48 & \textbf{5.99*} \\
 &  & TW4 &  & 7.29 & 6.38 & 8.11 & 6.39 & \textbf{5.70*} \\
 &  & TW8 &  & 6.67 & 5.92 & 7.38 & 5.97 & \textbf{5.29*} \\
     \cmidrule(lr){2-9}
 & \multirow[t]{6}{*}{UNI} & DL2 &  & 2.20 & 0.12 & 0.16 & 0.17 & \textbf{0.12} \\
 &  & DL4 &  & 4.65 & 2.63 & 3.16 & 2.63 & \textbf{2.55} \\
 &  & DL8 &  & 6.77 & 5.10 & 6.39 & 5.12 & \textbf{5.02} \\
     \cmidrule(lr){3-9}
 &  & TW2 &  & 9.95 & 9.91 & 11.15 & 10.25 & \textbf{9.21*} \\
 &  & TW4 &  & 9.98 & 10.14 & 11.14 & 10.75 & \textbf{9.53*} \\
 &  & TW8 &  & 9.63 & 10.15 & 10.64 & 10.77 & \textbf{9.26} \\
 \midrule
\multirow[t]{12}{*}{RC} & \multirow[t]{6}{*}{HOM} & DL2 &  & 3.33 & 2.29 & 2.76 & \textbf{2.08} & 2.23 \\
 &  & DL4 &  & 7.79 & 6.33 & 7.87 & \textbf{5.94} & 6.19 \\
 &  & DL8 &  & 10.83 & 9.76 & 12.27 & 9.39 & \textbf{9.36} \\
     \cmidrule(lr){3-9}
 &  & TW2 &  & 12.84 & 11.95 & 13.09 & 12.09 & \textbf{10.67*} \\
 &  & TW4 &  & 13.57 & 12.95 & 14.48 & 12.92 & \textbf{11.65*} \\
 &  & TW8 &  & 13.07 & 13.00 & 14.42 & 13.04 & \textbf{11.59*} \\
     \cmidrule(lr){2-9}
 & \multirow[t]{6}{*}{UNI} & DL2 &  & 4.41 & 2.49 & 2.75 & 2.37 & \textbf{2.34} \\
 &  & DL4 &  & 10.47 & 7.45 & 8.63 & 7.43 & \textbf{7.28} \\
 &  & DL8 &  & 14.37 & 11.60 & 13.34 & 11.59 & \textbf{11.25} \\
     \cmidrule(lr){3-9}
 &  & TW2 &  & 14.21 & 15.24 & 15.33 & 15.45 & \textbf{13.61*} \\
 &  & TW4 &  & \textbf{14.59} & 16.19 & 16.38 & 16.30 & 14.67 \\
 &  & TW8 &  & \textbf{14.37} & 16.10 & 16.38 & 16.54 & 14.72 \\
\midrule
Average & & & & 9.58 & 8.76 & 9.66 & 8.83 & \textbf{8.05} \\
\bottomrule
\end{tabular}
\end{table}

Several observations emerge from the benchmark results. 
First, under a single threshold consensus function (DSHH or ICD-postpone), it is generally better to operate with a dispatch threshold rather than a postponement threshold.
We believe that this is due to the fact that the postponement threshold only depends on sampled future request realizations.
In contrast, dispatch thresholds also take into account immediate costs that are based on already-known requests across all scenarios. Therefore, dispatch thresholds make for a more reliable consensus function with less noise.
Second, there is value in dispatching and postponing simultaneously: there is a notable difference in performance between iterative methods with one-sided decisions (DSHH and ICD-postpone) and iterative methods with two-sided decisions (ICD-Hamming and ICD-double). 
When we reduce the total expected number of requests from 600 to 450 and further to 300 (Table~\ref{tab:arrivals}), these differences persist.
This shows the benefits of simultaneous dispatching and postponement decisions, even when the scenario instances are considerably smaller.
For smaller instances with 300 expected arrivals, the non-parametric ICD-Hamming solution obtains very good results. 
As the instances increase in size, however, the single-scenario solution that ICD-Hamming uses is outperformed by the threshold-based ICD-double algorithm, which uses information from all scenario solutions.
The extended results for 300 and 450 total number of expected requests can be found in Appendix~\ref{app:results}.
Third, the ICD algorithms (except ICD-postpone) generally perform better than RH, showing that there is value in considering multiple scenarios and iterations.
Among instance classes with unimodal arrival processes and time windows, the non-iterative RH policy obtains good results.
In such instances, many requests can be postponed until the final epochs.
This may lead to very large scenario instances, which under tight time limits, may favor a single, optimized routing solution.
Nevertheless, ICD-double still achieves better results than RH in the majority of the instance classes (34 out of 36).

\begin{table}[ht]
\caption{
The average percentage gap to the hindsight solution for different numbers of total expected number of requests.
}
\label{tab:arrivals}

\centering
\begin{tabular}{lrrrrr}
\toprule
             & RH & DSHH & ICD-postpone & ICD-Hamming & ICD-double \\
             \midrule
300 arrivals & 12.58 & 10.57 & 11.84 & \textbf{9.93}  & 10.03 \\
450 arrivals & 11.12 & 9.61 & 10.76 & 9.36 & \textbf{8.96} \\
600 arrivals & 9.58  & 8.76  & 9.66  & 8.83  & \textbf{8.05}  \\
\bottomrule
\end{tabular}
\end{table}

Next, we conduct further experiments to shed light on the inner workings of ICD and highlight how important parameter settings affect the algorithm's performance.
We restrict these additional experiments to the benchmark instances with unimodal arrival processes and time windows.
All parameters are as in Section~\ref{sec:parameter-settings} unless mentioned otherwise. 

\begin{figure}[ht]
    \centering
    \includegraphics[width=0.8\linewidth]{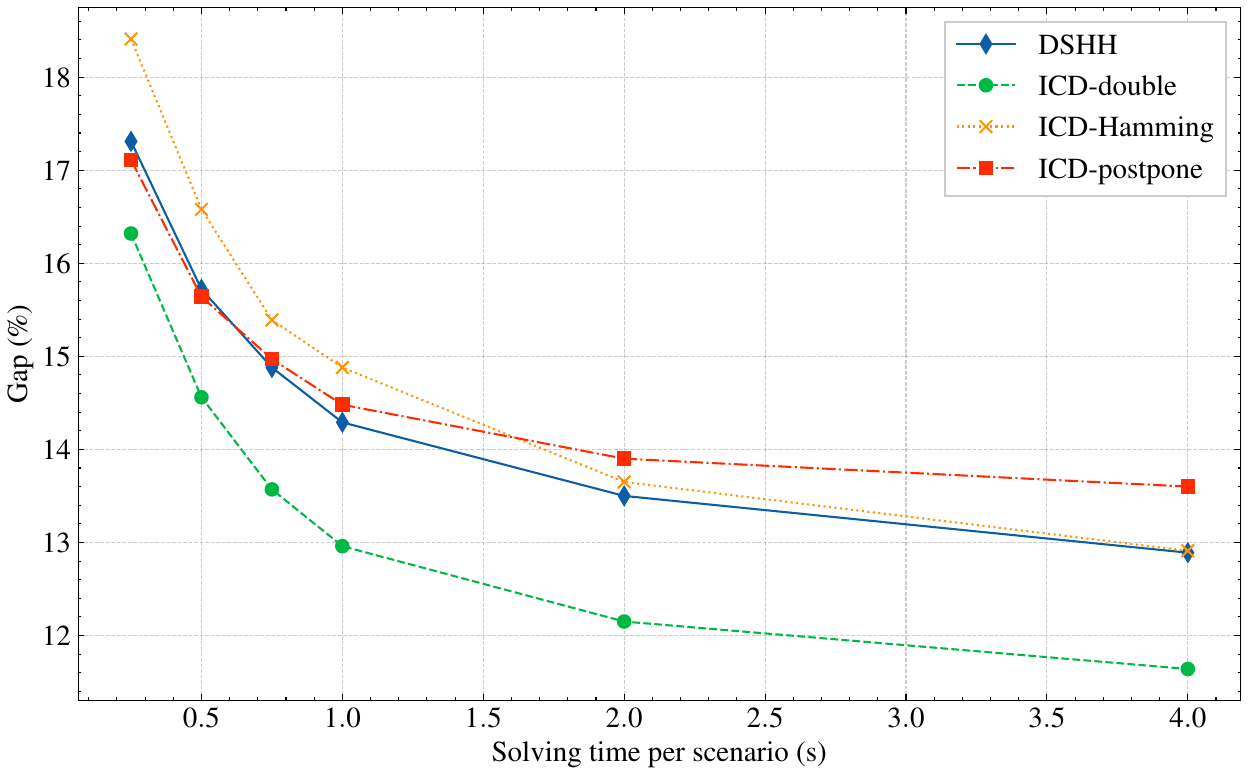}
    \caption{
    Average percentage gaps to the hindsight solution versus the solving time in seconds per scenario.
    }
    \label{fig:icd-runtime}
\end{figure}

\paragraph{Solving time per scenario} 
Varying the time allocated to solving sample scenarios, we evaluate the impact of ICD on the solution of sample scenarios and the impact on the overall solution of the complete dynamic problem. 
We vary the amount of time spent solving a single scenario between 0.25, 0.5, 0.75, 1, 2, and 4 seconds (recall that the default is 1 second).
The time to compute the direct cost remains fixed at 30 seconds for all considered scenario time limits.
We compare all results to the hindsight solutions obtained after 600 seconds.
The results are presented in Figure~\ref{fig:icd-runtime}, which shows the average gap with respect to the hindsight solution versus the solving time for a single scenario.
As the sample scenarios are solved with more time, the resulting solutions provide a better trade-off in minimizing the immediate and future costs.
When increasing the solving time from 0.25 seconds to 1 second, the average percentage gap decreases between 2.5-4\% across all methods. Increasing the solving time even further continues to improve the solutions, but at diminishing returns. These observations are consistent with the convergence curve of our static solver on general VRP instances, showing that the improvement from obtaining better static solutions translates to its dynamic counterpart.

\begin{figure}[ht!]
    \centering
    \includegraphics[width=0.8\linewidth]{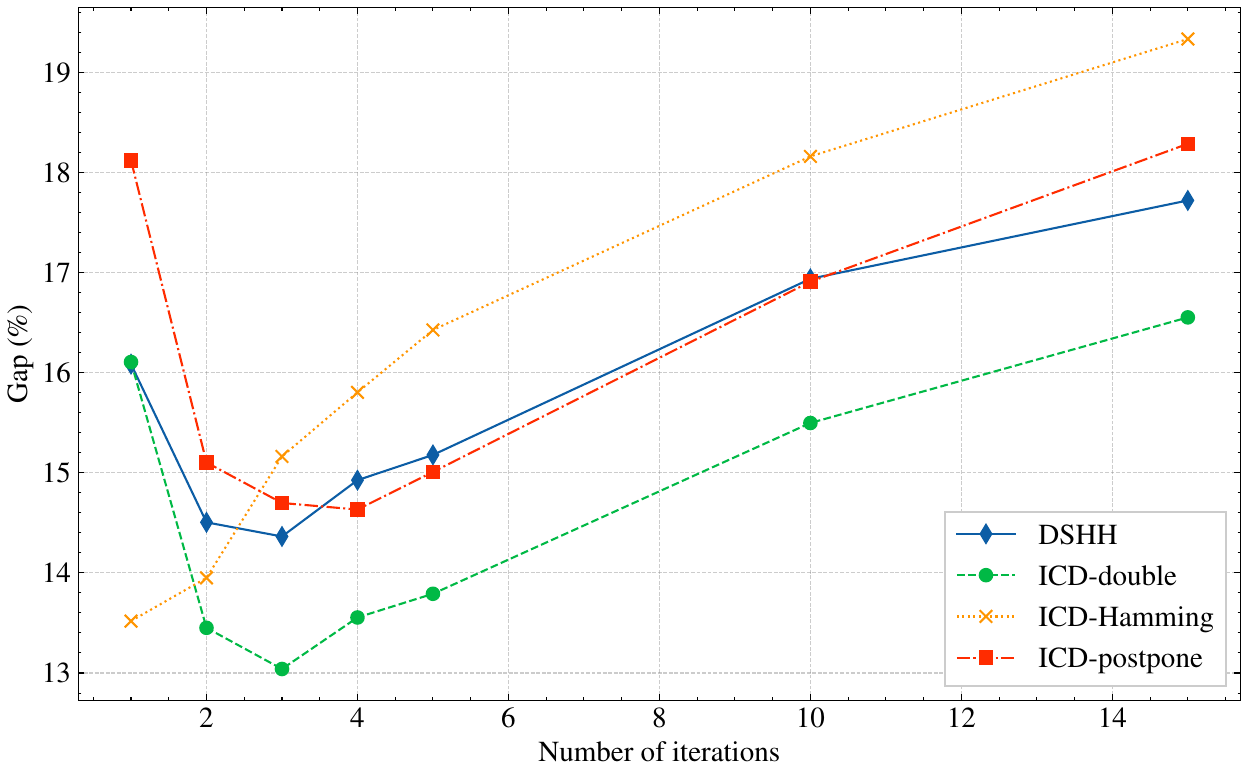}
    \caption{
    Average percentage gaps to the hindsight solution versus the number of iterations in the ICD method.
    }
    \label{fig:icd-iterations}
\end{figure}

\paragraph{Number of iterations.} 
In the next experiment, we alter the number of iterations while keeping all other ICD and consensus function parameters constant.
We test the number of iterations using the values 1, 2, 3, 5, 10, and 15. 
Figure~\ref{fig:icd-iterations} illustrates the average gap with respect to the hindsight solution versus the number of iterations.
For the threshold-based consensus functions, the optimal number of iterations is found to be 3. 
This also supports the choices made in the initial parameter value selection of Section~\ref{sec:parameter-settings}.
For ICD-Hamming, the ideal number of iterations in this experiment is one. 
We note that on the complete set of instances, as presented in Table~\ref{tab:main-results}, the overall performance of ICD-Hamming with three iterations was better than with one iteration, due to its superior performance on the deadline instances. 
Remarkably, a single iteration for the ICD-Hamming method performs substantially better than the other methods using only one iteration.
This can be in part explained by the fact that the Hamming distance consensus function selects a complete scenario solution for which the immediate reward is known.
The threshold-based methods select requests that do not necessarily form a coherent route plan, as the requests are selected based on their relative dispatching frequency, not how often they are dispatched together. 
The dispatching action after one iteration thus does not necessarily result in efficient route plans.
Instead, the threshold-based methods refine this dispatching action in subsequent iterations, given that part of the dispatching action is already determined.
As the number of iterations increases beyond three, the trade-off between refining the dispatching action and solving scenarios becomes unfavorable.
For example, using 10 iterations with the selected parameters results in solving 300 scenarios instead of 90 scenarios within the same time budget.
This aligns with the earlier observation that allocating too little time to solving scenarios leads to worse overall performance.

\subsection{Results on the \textit{EURO meets NeurIPS 2022} vehicle routing competition instances}
\label{sec:euro-neurips}
This final section evaluates the ICD-double algorithm on instances from the \textit{EURO meets NeurIPS 2022} vehicle routing competition \citep{Kool2022a}.
In the following, we abbreviate the competition name by \textit{EURO-NeurIPS}.

The EURO-NeurIPS competition presented two VRP variants that needed to be solved: a static VRPTW and a dynamic version thereof.
We focus on the dynamic problem variant, which is also precisely what we formulated as the DDWP in this work.
The main differences between the instances used in the competition versus our work are the characteristics of the static instances and how the requests are sampled.
In the competition, an instance of the DDWP is created by sampling at each decision epoch 100 requests from a static VRPTW instance.
However, sampled requests that cannot be served on time are disregarded, meaning that the number of requests to be served is often much lower, especially in later epochs.
Depending on the static instance, the total number of epochs ranges between 5 and 9, leading to a total number of requests ranging between 350 to 600.

After a four-month qualification phase, the ten best algorithms were evaluated on 100 final hidden instances with a time limit of 120 seconds per epoch.
We evaluate ICD-double on these 100 instances for the same time limits of 120 seconds, as well as longer time limits of 180 and 600 seconds, and compare with results from the competition.
In the competition, the algorithms were evaluated on a single core of an Intel Xeon Gold 5118 processor with a CPU PassMark score of 1808, which is slightly lower than the CPU PassMark score of 2014 of the CPU that we used. 
To provide a fair comparison, we compensate for the CPU differences by reducing our time limits by a factor of $\frac{1808}{2014}$.
Moreover, we did not change our parameters or implementation from the ones used in the benchmark instances in previous sections.

\begin{figure}
    \centering
    \includegraphics[width=0.9\linewidth]{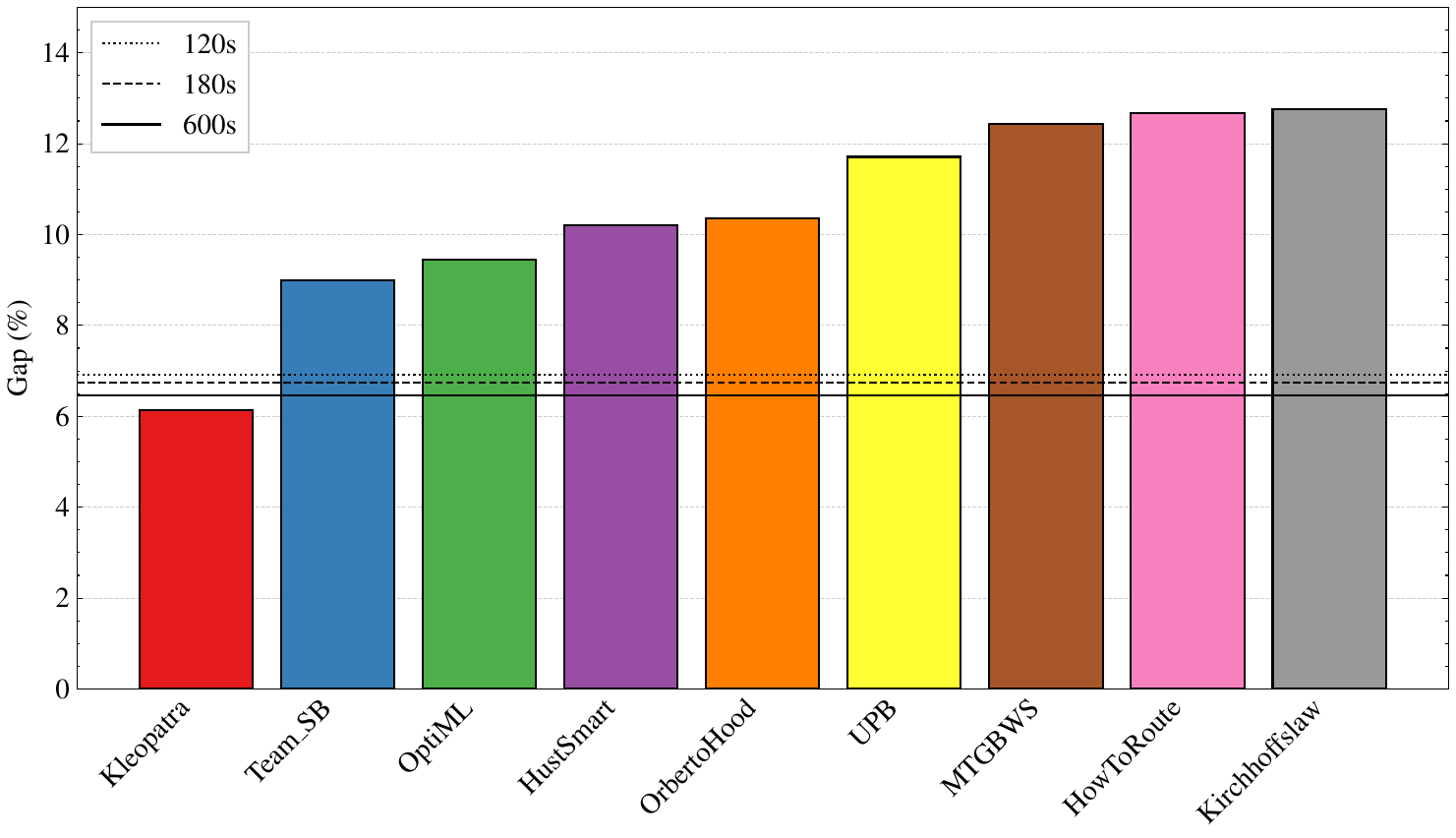}
    \caption{
    Results for the final 100 instances of the EURO-NeurIPS competition.
    The lines indicate the average percentage gap to the hindsight solution obtained by ICD-double for a given time limit (120, 180, or 600 seconds).
    }
    \label{fig:euro-neurips}
\end{figure}

Figure~\ref{fig:euro-neurips} shows the average percentage gap to the hindsight solution for each final team, along with three lines of the ICD-double algorithm for a given epoch time limit.
The ICD-double algorithm with the 120-second time limit would have ranked second and has a 0.8\% gap over the best result.
With a slightly longer time limit of 180 seconds, ICD-double is only 0.6\% away, and increasing the time limit even further to 600 seconds decreases this gap to 0.3\%.
These results align with our earlier observations (cf. Figure~\ref{fig:icd-runtime}), where increasing the time to solve a single scenario improves overall performance.

We briefly summarize the contribution of the top six teams.
It is important to note that all teams were provided with a state-of-the-art VRPTW solver by \cite{Kool2022b}. 
PyVRP is built on top of the same static solver with marginal performance differences, meaning that differences between our ICD-double and other methods are mainly due to differences in the dynamic aspect of the algorithm.
Kleopatra~\citep{Baty2024} proposed a machine learning pipeline to predict a prize associated with each request, quantifying how favorable it is to dispatch a request in the current epoch.
They trained their model on the solutions of hindsight instances and solved a prize-collecting variant of the VRPTW in each decision epoch using the predicted prizes.
Team\_SB~\citep{kwon2022vehicle} also proposed a machine learning-based approach. 
They trained two neural networks, one that predicts priorities for requests and one that selects requests that must be dispatched.
With the predicted priority values, they solve a VRPTW with a modified cost function from which they determine which requests to dispatch.
We, the author team, participated under the team name OptiML, and our submitted approach was a variant of the ICD-postpone algorithm with parameters tuned for the EURO-NeurIPS instances that were known to us during the qualification phase of the competition.
HustSmart~\citep{li2022vehicle} heuristically defined penalties for each request based on the static instance data. 
They then solved a variant of the VRP with penalties and dropping visits, which can be seen as a variant of the prize-collecting VRP.
OrbertoHood~\citep{hildebrandt2022vehicle} was the only other finalist who also used an online sampling-based approach. 
Their approach samples a set of scenarios and then solves a sample average approximation problem to estimate the marginal costs of dispatching or postponing each request. 
They solved the problem using column generation and used the resulting solution to determine which requests to dispatch.
Finally, UPB~\citep{Dieter2023} proposed a regret policy that analytically computes a regret value for each optional request.
The regret value is the potential value that would be missed if a request is dispatched in the current epoch, and their algorithm postpones all requests with a regret value above a certain threshold.

\section{Conclusions}\label{sec:conclusion}
This paper studies a variant of the dynamic dispatch waves problem (DDWP) in which the goal is to minimize the total routing costs while serving all incoming delivery requests on time.
We propose the iterative conditional dispatch (ICD) algorithm, which iteratively solves sample scenarios to classify requests to be dispatched, postponed, or leaves this decision to be determined for the next iterations.
The most competitive variant of the algorithm uses a double threshold consensus function, which classifies requests based on their dispatching frequency in the resolved sample scenarios.
We additionally develop a parameter-free consensus function based on the Hamming distance.
Numerical experiments on a wide variety of instances demonstrate that ICD with the Hamming distance consensus function generates competitive solutions compared to the dynamic stochastic hedging heuristic of~\cite{Hvattum2006}, while the ICD with double threshold improves over both methods by as much as 2\% on average for some instance classes.
Additional experiments show that spending sufficient computational time to solve scenario instances is important to obtain overall good performance, but that there are clear diminishing returns.
We find that the threshold-based variants of ICD obtain their best performance with three iterations, while increasing the number of iterations typically results in inferior solutions due to the reduced time available for solving scenarios.
Finally, we show using instances from the \textit{EURO meets NeurIPS 2022} vehicle routing competition that our ICD method with double thresholds achieves excellent results.

The performance and simplicity of the ICD algorithm show a promising avenue for extension to other dynamic vehicle routing problem variants.
In future work, we intend to apply ICD to other problems, including routing problems such as the same-day delivery problem with continuous arrivals~\citep{Voccia2019} and dispatching problems with ad-hoc fleet arrivals \citep{arslan2019crowdsourced}.

\ACKNOWLEDGMENT{%
Leon Lan and Jasper van Doorn would like to thank TKI Dinalog, the Topsector Logistics, and the Dutch Ministry of Economic Affairs and Climate Policy for funding this project.
}

\clearpage
\begin{APPENDICES}
    \section{Additional results}
\label{app:results}
This section presents additional results to the numerical experiments of Section~\ref{sec:experiments}.
Table~\ref{tab:main-results-p-values} presents the p-values of the pair-wise T-tests for the benchmark results (Table~\ref{tab:main-results}).
Table~\ref{tab:main-results-300} and~\ref{tab:main-results-450} show the results of the benchmark instances for 300 and 450 expected number of arrivals, respectively.

\begin{table}[ht]
\caption{
    Overview of p-values of the pair-wise T-tests. The method with the best average gap in each row is used as the benchmark and represented with a dash (-).}
\label{tab:main-results-p-values}
\small
\centering
\begin{tabular}{llllrrrrr}
\toprule
\multicolumn{3}{c}{Instance} && \multicolumn{5}{c}{Methods} \\
\cmidrule{1-3} \cmidrule{5-9}
Topology & Arrival & TW &&  RH &  DSHH &  ICD-postpone &  ICD-Hamming &  ICD-double \\
\midrule
\multirow[t]{12}{*}{R} & \multirow[t]{6}{*}{HOM} & DL2 &  & \textless 0.001 & 0.019 & \textless 0.001 & - & 0.165 \\
 &  & DL4 &  & \textless 0.001 & \textless 0.001 & \textless 0.001 & 0.826 & - \\
 &  & DL8 &  & \textless 0.001 & 0.005 & \textless 0.001 & 0.445 & - \\
     \cmidrule(lr){3-9}
 &  & TW2 &  & \textless 0.001 & \textless 0.001 & \textless 0.001 & \textless 0.001 & - \\
 &  & TW4 &  & \textless 0.001 & \textless 0.001 & \textless 0.001 & \textless 0.001 & - \\
 &  & TW8 &  & \textless 0.001 & \textless 0.001 & \textless 0.001 & \textless 0.001 & - \\
\cmidrule(lr){2-9}
 & \multirow[t]{6}{*}{UNI} & DL2 &  & \textless 0.001 & 0.653 & \textless 0.001 & 0.348 & - \\
 &  & DL4 &  & \textless 0.001 & 0.002 & \textless 0.001 & 0.128 & - \\
 &  & DL8 &  & \textless 0.001 & \textless 0.001 & \textless 0.001 & \textless 0.001 & - \\
     \cmidrule(lr){3-9}
 &  & TW2 &  & \textless 0.001 & \textless 0.001 & \textless 0.001 & \textless 0.001 & - \\
 &  & TW4 &  & 0.119 & \textless 0.001 & \textless 0.001 & \textless 0.001 & - \\
 &  & TW8 &  & 0.726 & \textless 0.001 & \textless 0.001 & \textless 0.001 & - \\
\midrule
\multirow[t]{12}{*}{C} & \multirow[t]{6}{*}{HOM} & DL2 &  & \textless 0.001 & - & 0.014 & 0.044 & 0.911 \\
 &  & DL4 &  & \textless 0.001 & 0.132 & \textless 0.001 & - & 0.260 \\
 &  & DL8 &  & \textless 0.001 & 0.019 & \textless 0.001 & - & 0.075 \\
     \cmidrule(lr){3-9}
 &  & TW2 &  & \textless 0.001 & 0.001 & \textless 0.001 & \textless 0.001 & - \\
 &  & TW4 &  & \textless 0.001 & \textless 0.001 & \textless 0.001 & \textless 0.001 & - \\
 &  & TW8 &  & \textless 0.001 & \textless 0.001 & \textless 0.001 & \textless 0.001 & - \\
\cmidrule(lr){2-9}
 & \multirow[t]{6}{*}{UNI} & DL2 &  & \textless 0.001 & 0.936 & 0.006 & 0.006 & - \\
 &  & DL4 &  & \textless 0.001 & 0.452 & \textless 0.001 & 0.497 & - \\
 &  & DL8 &  & \textless 0.001 & 0.526 & \textless 0.001 & 0.429 & - \\
     \cmidrule(lr){3-9}
 &  & TW2 &  & \textless 0.001 & \textless 0.001 & \textless 0.001 & \textless 0.001 & - \\
 &  & TW4 &  & 0.002 & \textless 0.001 & \textless 0.001 & \textless 0.001 & - \\
 &  & TW8 &  & 0.017 & \textless 0.001 & \textless 0.001 & \textless 0.001 & - \\
\midrule
\multirow[t]{12}{*}{RC} & \multirow[t]{6}{*}{HOM} & DL2 &  & \textless 0.001 & 0.052 & \textless 0.001 & - & 0.170 \\
 &  & DL4 &  & \textless 0.001 & 0.022 & \textless 0.001 & - & 0.130 \\
 &  & DL8 &  & \textless 0.001 & 0.043 & \textless 0.001 & 0.849 & - \\
     \cmidrule(lr){3-9}
 &  & TW2 &  & \textless 0.001 & \textless 0.001 & \textless 0.001 & \textless 0.001 & - \\
 &  & TW4 &  & \textless 0.001 & \textless 0.001 & \textless 0.001 & \textless 0.001 & - \\
 &  & TW8 &  & \textless 0.001 & \textless 0.001 & \textless 0.001 & \textless 0.001 & - \\
\cmidrule(lr){2-9}
 & \multirow[t]{6}{*}{UNI} & DL2 &  & \textless 0.001 & 0.169 & \textless 0.001 & 0.826 & - \\
 &  & DL4 &  & \textless 0.001 & 0.317 & \textless 0.001 & 0.359 & - \\
 &  & DL8 &  & \textless 0.001 & 0.092 & \textless 0.001 & 0.109 & - \\
     \cmidrule(lr){3-9}
 &  & TW2 &  & 0.004 & \textless 0.001 & \textless 0.001 & \textless 0.001 & - \\
 &  & TW4 &  & - & \textless 0.001 & \textless 0.001 & \textless 0.001 & 0.736 \\
 &  & TW8 &  & - & \textless 0.001 & \textless 0.001 & \textless 0.001 & 0.159 \\
\bottomrule
\end{tabular}
\end{table}

\begin{table}[ht]
\caption{
    Overview of results with 300 expected total number of requests.
    The values represent the average percentage gap with respect to the hindsight solution over 200 different instances.
    The best average gap in each row is marked in bold.
    An asterisk (*) indicates the method is statistically significantly better than the other methods for that instance at the 0.05 level (with Bonferroni correction).}
\label{tab:main-results-300}
\small
\centering
\begin{tabular}{llllrrrrr}
\toprule
\multicolumn{3}{c}{Instance} && \multicolumn{5}{c}{Methods} \\
\cmidrule{1-3} \cmidrule{5-9}
Topology & Arrival & TW &&  RH &  DSHH &  ICD-postpone &  ICD-Hamming &  ICD-double \\
\midrule
\multirow[t]{12}{*}{R} & \multirow[t]{6}{*}{HOM} & DL2 &  & 6.23 & 4.00 & 4.23 & \textbf{3.95} & 4.01 \\
 &  & DL4 &  & 11.58 & 8.73 & 9.94 & \textbf{8.07*} & 8.65 \\
 &  & DL8 &  & 15.25 & 12.81 & 15.30 & \textbf{11.46*} & 12.54 \\
     \cmidrule(lr){3-9}
 &  & TW2 &  & 17.05 & 14.31 & 16.23 & 13.61 & \textbf{13.32} \\
 &  & TW4 &  & 17.65 & 15.20 & 17.25 & 13.98 & \textbf{13.95} \\
 &  & TW8 &  & 16.70 & 14.80 & 16.95 & 13.82 & \textbf{13.59} \\
     \cmidrule(lr){2-9}
 & \multirow[t]{6}{*}{UNI} & DL2 &  & 6.50 & 4.21 & 4.49 & \textbf{3.98} & 4.15 \\
 &  & DL4 &  & 12.35 & 9.53 & 10.69 & \textbf{8.84} & 9.24 \\
 &  & DL8 &  & 16.26 & 13.85 & 15.81 & \textbf{12.81} & 13.36 \\
     \cmidrule(lr){3-9}
 &  & TW2 &  & 16.98 & 15.61 & 16.62 & 14.97 & \textbf{14.53} \\
 &  & TW4 &  & 17.05 & 16.63 & 17.40 & 15.70 & \textbf{15.33} \\
 &  & TW8 &  & 16.56 & 16.93 & 16.81 & 15.81 & \textbf{15.35} \\
 \midrule
\multirow[t]{12}{*}{C} & \multirow[t]{6}{*}{HOM} & DL2 &  & 3.18 & 1.41 & 1.58 & 1.43 & \textbf{1.41} \\
 &  & DL4 &  & 6.96 & 4.18 & 4.63 & \textbf{3.99} & 4.08 \\
 &  & DL8 &  & 8.49 & 6.11 & 7.19 & \textbf{5.61*} & 6.08 \\
     \cmidrule(lr){3-9}
 &  & TW2 &  & 11.28 & 8.91 & 10.73 & \textbf{8.39} & 8.53 \\
 &  & TW4 &  & 10.63 & 8.59 & 10.57 & \textbf{7.91} & 8.09 \\
 &  & TW8 &  & 9.79 & 8.11 & 9.85 & \textbf{7.51} & 7.64 \\
     \cmidrule(lr){2-9}
 & \multirow[t]{6}{*}{UNI} & DL2 &  & 4.20 & 1.49 & 1.51 & \textbf{1.49} & 1.49 \\
 &  & DL4 &  & 8.70 & 5.61 & 6.05 & \textbf{5.29} & 5.38 \\
 &  & DL8 &  & 10.94 & 8.26 & 9.61 & \textbf{7.72} & 8.13 \\
     \cmidrule(lr){3-9}
 &  & TW2 &  & 13.25 & 11.92 & 13.52 & \textbf{11.17} & 11.26 \\
 &  & TW4 &  & 13.06 & 11.93 & 13.51 & 11.52 & \textbf{11.13} \\
 &  & TW8 &  & 12.77 & 11.80 & 12.66 & 11.39 & \textbf{11.00} \\
 \midrule
\multirow[t]{12}{*}{RC} & \multirow[t]{6}{*}{HOM} & DL2 &  & 5.96 & 3.63 & 4.03 & \textbf{3.56} & 3.73 \\
 &  & DL4 &  & 11.48 & 8.33 & 9.53 & \textbf{7.67*} & 8.28 \\
 &  & DL8 &  & 15.36 & 12.60 & 15.19 & \textbf{11.65*} & 12.55 \\
     \cmidrule(lr){3-9}
 &  & TW2 &  & 17.29 & 14.43 & 16.52 & 13.76 & \textbf{13.56} \\
 &  & TW4 &  & 17.77 & 14.67 & 17.09 & 14.14 & \textbf{13.85} \\
 &  & TW8 &  & 17.12 & 14.48 & 17.04 & 13.61 & \textbf{13.52} \\
     \cmidrule(lr){2-9}
 & \multirow[t]{6}{*}{UNI} & DL2 &  & 6.15 & 4.21 & 4.41 & \textbf{3.96} & 4.24 \\
 &  & DL4 &  & 11.83 & 9.56 & 10.69 & \textbf{8.78*} & 9.47 \\
 &  & DL8 &  & 15.28 & 13.57 & 15.64 & \textbf{12.73} & 13.27 \\
     \cmidrule(lr){3-9}
 &  & TW2 &  & 17.10 & 16.21 & 17.41 & 15.26 & \textbf{15.01} \\
 &  & TW4 &  & 17.32 & 16.72 & 18.00 & 15.99 & \textbf{15.58} \\
 &  & TW8 &  & 16.72 & 17.19 & 17.69 & 15.98 & \textbf{15.80} \\
\midrule
Average & & & & 12.58 & 10.57 & 11.84 & \textbf{9.93} & 10.03 \\
\bottomrule
\end{tabular}
\end{table}

\begin{table}[ht]
\caption{
    Overview of results with 450 expected total number of requests.
    The values represent the average relative gap (in percentages) with respect to the hindsight solution over 200 different instances.
    The best average gap in each row is marked in bold.
    An asterisk (*) indicates the method is statistically significantly better than the other methods for that instance at the 0.05 level (with Bonferroni correction).}
\label{tab:main-results-450}
\small
\centering
\begin{tabular}{llllrrrrr}
\toprule
\multicolumn{3}{c}{Instance} && \multicolumn{5}{c}{Methods} \\
\cmidrule{1-3} \cmidrule{5-9}
Topology & Arrival & TW &&  RH &  DSHH &  ICD-postpone &  ICD-Hamming &  ICD-double \\
\midrule
\multirow[t]{12}{*}{R} & \multirow[t]{6}{*}{HOM} & DL2 &  & 5.04 & 3.36 & 3.80 & \textbf{3.18} & 3.34 \\
 &  & DL4 &  & 9.76 & 7.76 & 9.17 & \textbf{7.12*} & 7.53 \\
 &  & DL8 &  & 13.27 & 11.31 & 13.72 & \textbf{10.58} & 10.87 \\
     \cmidrule(lr){3-9}
 &  & TW2 &  & 14.40 & 13.05 & 15.02 & 12.95 & \textbf{12.04*} \\
 &  & TW4 &  & 14.68 & 14.15 & 15.84 & 13.79 & \textbf{12.86*} \\
 &  & TW8 &  & 14.15 & 14.06 & 15.84 & 14.05 & \textbf{12.47*} \\
     \cmidrule(lr){2-9}
 & \multirow[t]{6}{*}{UNI} & DL2 &  & 5.30 & 3.54 & 3.83 & \textbf{3.48} & 3.49 \\
 &  & DL4 &  & 10.97 & 8.67 & 9.74 & 8.34 & \textbf{8.30} \\
 &  & DL8 &  & 14.48 & 12.89 & 14.56 & 12.60 & \textbf{12.36} \\
     \cmidrule(lr){3-9}
 &  & TW2 &  & 15.55 & 15.89 & 16.00 & 15.86 & \textbf{14.28*} \\
 &  & TW4 &  & 15.79 & 16.84 & 16.98 & 16.48 & \textbf{15.18*} \\
 &  & TW8 &  & \textbf{14.85} & 16.85 & 16.54 & 16.65 & 15.02 \\
 \midrule
\multirow[t]{12}{*}{C} & \multirow[t]{6}{*}{HOM} & DL2 &  & 1.77 & 0.35 & 0.37 & 0.34 & \textbf{0.33} \\
 &  & DL4 &  & 4.75 & 2.67 & 3.31 & \textbf{2.50} & 2.74 \\
 &  & DL8 &  & 6.19 & 4.57 & 5.83 & \textbf{4.08*} & 4.48 \\
     \cmidrule(lr){3-9}
 &  & TW2 &  & 9.24 & 7.78 & 9.55 & 7.38 & \textbf{7.31} \\
 &  & TW4 &  & 8.78 & 7.51 & 9.40 & 7.19 & \textbf{7.01} \\
 &  & TW8 &  & 8.25 & 6.99 & 8.68 & 6.73 & \textbf{6.31*} \\
     \cmidrule(lr){2-9}
 & \multirow[t]{6}{*}{UNI} & DL2 &  & 2.91 & 0.40 & 0.41 & \textbf{0.37} & 0.38 \\
 &  & DL4 &  & 6.41 & 4.17 & 4.69 & \textbf{4.01} & 4.06 \\
 &  & DL8 &  & 8.51 & 6.50 & 7.97 & \textbf{6.38} & 6.41 \\
     \cmidrule(lr){3-9}
 &  & TW2 &  & 11.60 & 10.75 & 12.30 & 10.70 & \textbf{10.10*} \\
 &  & TW4 &  & 11.59 & 11.00 & 12.28 & 11.08 & \textbf{10.32*} \\
 &  & TW8 &  & 11.27 & 10.84 & 11.40 & 10.87 & \textbf{9.87*} \\
 \midrule
\multirow[t]{12}{*}{RC} & \multirow[t]{6}{*}{HOM} & DL2 &  & 4.70 & 2.85 & 3.36 & \textbf{2.66} & 2.80 \\
 &  & DL4 &  & 9.15 & 7.19 & 8.84 & \textbf{6.52} & 6.95 \\
 &  & DL8 &  & 12.64 & 10.79 & 13.54 & \textbf{10.03} & 10.39 \\
     \cmidrule(lr){3-9}
 &  & TW2 &  & 14.95 & 13.07 & 15.06 & 12.86 & \textbf{11.93*} \\
 &  & TW4 &  & 15.64 & 13.69 & 15.81 & 13.47 & \textbf{12.57*} \\
 &  & TW8 &  & 15.13 & 13.61 & 15.64 & 13.51 & \textbf{12.48*} \\
     \cmidrule(lr){2-9}
 & \multirow[t]{6}{*}{UNI} & DL2 &  & 5.68 & 2.96 & 3.29 & \textbf{2.86} & 2.95 \\
 &  & DL4 &  & 12.05 & 8.12 & 9.48 & \textbf{7.91} & 8.01 \\
 &  & DL8 &  & 15.68 & 12.27 & 14.18 & \textbf{11.68} & 11.79 \\
     \cmidrule(lr){3-9}
 &  & TW2 &  & 18.24 & 15.82 & 16.64 & 15.63 & \textbf{14.64*} \\
 &  & TW4 &  & 18.52 & 16.87 & 17.32 & 16.48 & \textbf{15.60*} \\
 &  & TW8 &  & 18.38 & 16.73 & 17.06 & 16.48 & \textbf{15.30*} \\
\midrule
Average & & & & 11.12 & 9.61 & 10.76 & 9.36 & \textbf{8.96} \\
\bottomrule
\end{tabular}
\end{table}

\clearpage
\section{Implementation details of solving static VRPTW-DW}
\label{app:implementation-details}
To solve static VRPTW-DWs, we rely on PyVRP \citep{wouda_pyvrp_2023}, which implements a hybrid genetic search (HGS) algorithm that builds upon HGS-CVRP \citep{vidal2022hgscvrp}.
We first present a general outline of HGS in Section~\ref{sec:hgs-algorithm-outline}, and we then describe how to extend it to solve routing problems with dispatch windows in Section~\ref{sec:hgs-dispatch-windows}.
Section~\ref{sec:hgs-parameter-tuning} presents our parameter tuning procedure and the resulting parameters.

\subsection{Algorithm outline}\label{sec:hgs-algorithm-outline}
HGS combines elements of a genetic algorithm and a local search algorithm, resulting in an effective synergy between the genetic algorithm's ability to explore a diverse set of solutions and local search to identify local optima.
In particular, HGS manages a population of candidate solutions, which are iteratively recombined and improved using a local search procedure.
An outline of the algorithm is given in Algorithm~\ref{alg:hgs-pyvrp}.

\begin{algorithm}
\caption{PyVRP's hybrid genetic search}
\label{alg:hgs-pyvrp}
\begin{algorithmic}[1]
    \State \textbf{Input}: initial solutions $s_1, \dots, s_p$
    \State Initialize $s^*$ to the initial solution with the best objective value
    \While{time limit is not met}
        \State Select two parent solutions from the population using binary tournament
        \State Apply a crossover operator on the parent solutions to generate an offspring solution
        \State Improve the offspring using a local search algorithm to obtain a candidate solution $s^c$
        \State Add the candidate solution to the population
        \If{$s^c$ has better objective value than $s^*$}
            \State $s^* \leftarrow s^c$
        \EndIf
        \If{population size exceeds maximum size}
            \State Remove the solutions with the lowest fitness until the population is at minimum size
        \EndIf
    \EndWhile
    \State \textbf{Output}: the best-found solution $s^*$
\end{algorithmic}
\end{algorithm}

The HGS algorithm works as follows. HGS initializes its population by a set of solutions given as input to the algorithm.
In each iteration, the algorithm selects two existing parent solutions from the population using a binary tournament, favoring solutions with higher fitness (i.e., a weighted combination of cost and diversity value).
A crossover operator then combines the two parent solutions to generate an offspring solution that inherits features from both parents.
After the crossover, the offspring solution is further improved using a local search procedure.
The resulting candidate solution is first added to the population.
If the candidate solution improves over the best solution found so far, it is also registered as the new best solution.
Upon reaching the maximum population size, a survivor selection mechanism removes the least fit solutions until the population is back at the minimum size.
The algorithm continues until a provided stopping criterion is met, returning the best solution found.

\subsection{Dispatch windows extension}\label{sec:hgs-dispatch-windows}
In this section, we describe how to extend PyVRP to solve VRPs with dispatch windows.
PyVRP implements the path concatenation schemes of \cite{vidal2013hgs}, which enables it to evaluate local search moves involving time windows in amortized $O(1)$ time.
The scheme relies on infeasible solutions that incur additional penalties based on \textit{time-warping} \citep{Nagata2010}: late arrivals at customer locations are warped back in time and penalized.

We first outline the concatenation scheme as introduced in \cite{vidal2013hgs} and later we describe how we extend this scheme to take into account dispatch windows.

Local search moves can be seen as a separation of routes into subsequences of location visits, which are then concatenated together into new routes.
Define $\sigma = (\sigma_{0}, \dots, \sigma_{l})$ as a sequence of visits (including depot).
The operator $\oplus$ concatenates two subsequences $\sigma, \sigmap$ together, forming a new sequence of visits.
For each subsequence $\sigma$, we maintain the following data: the minimum duration $D(\sigma)$, the time warp $TW(\sigma)$, the earliest $E(\sigma)$ and latest $L(\sigma)$ departure time to the first location with minimum duration and minimum time warp use, the cumulative distance $C(\sigma)$ and the cumulative load $Q(\sigma)$. For a sequence $\sigma^0 = \{i \}$ consisting of a single visit $i \in V$, the data is initialized as $D(\sigma^0) = \mu_i, TW(\sigma^0) = 0, E(\sigma^0) = e_i, L(\sigma^0) = l_i, C(\sigma^0) = 0$, and $Q(\sigma^0) = q_i$.

The following proposition shows how to compute the same data for a concatenation of sequences:
\begin{proposition}[Concatenation of two sequences \cite{vidal2013hgs}]
\label{prop:concat-scheme}
Let $\sigma = (\sigma_i, \dots, \sigma_j)$ and $\sigmap = (\sigmap_{i^\prime}, \dots, \sigmap_{j^{\prime}})$ be two subsequences of visits. The concatenated subsequence $\sigma \oplus \sigmap$ is characterized by the following data:
\begin{align}
& D\left(\sigma \oplus \sigmap\right)=D(\sigma)+D\left(\sigmap\right)+\delta_{\sigma_j \sigma_{i^\prime}}+\Delta_{W T} \\
& TW\left(\sigma \oplus \sigmap\right)=TW(\sigma)+ TW\left(\sigmap\right)+\Delta_{TW} \\
& E\left(\sigma \oplus \sigmap\right)=\max \left\{E\left(\sigmap\right)-\Delta, E(\sigma)\right\}-\Delta_{W T} \\
& L\left(\sigma \oplus \sigmap\right)=\min \left\{L\left(\sigmap\right)-\Delta, L(\sigma)\right\}+\Delta_{T W} \\
& C\left(\sigma \oplus \sigmap\right)=C(\sigma)+C\left(\sigmap\right)+c_{\sigma_j \sigmap_{i^\prime}} \\
& Q\left(\sigma \oplus \sigmap\right)=Q(\sigma)+Q\left(\sigmap\right)
\end{align}
where $\Delta=D(\sigma)-T W(\sigma)+\delta_{\sigma_j \sigmap_{i^\prime}}, \Delta_{W T}=\max \left\{E\left(\sigmap\right)-\Delta-L(\sigma), 0\right\}$ and
$\Delta_{T W}=\max \{E(\sigma)+\Delta-L(\sigmap), 0\}$.
\end{proposition}

Local search moves, which generally consist of up to five concatenation procedures, can be evaluated in constant time by computing the relevant data using this proposition.

We now describe how the concatenation scheme can be extended for dispatch windows.
For a sequence $\sigma$, we additionally compute the earliest $R^-(\sigma)$ and the latest $R^+(\sigma)$ dispatch time.
For a sequence $\sigma^0 = \{i\}$ consisting of a single visit $i \in V$, the data is initialized as $R^-(\sigma^0)=r^-_i$ and $R^+(\sigma^0) = r^+$.
The concatenation scheme of two subsequences $\sigma$ and $\sigmap$ is extended using the following expressions:
\begin{align}
& R^-\left(\sigma \oplus \sigmap\right) = \max \{R^- (\sigma), R^- (\sigmap)\} \\
& R^+\left(\sigma \oplus \sigmap\right) = \min \{R^+ (\sigma), R^+ (\sigmap)\}
\end{align}

Having introduced the dispatch window data, we aim to introduce additional time warp penalties when there is a violation of the time window data $E(\sigma), L(\sigma)$ and dispatch window data $R^-(\sigma), R^+(\sigma)$.
To this end, we rely on the following proposition:
\begin{proposition}[\cite{vidal2013hgs}]
\label{prop:starting-time}
    There exists a set of starting dates $t$ that minimize both the time warp use and the duration to service a sequence $\sigma$. This set is a segment, notated as $\mathcal{T}^{\text{MIN}}(\sigma) = \left[ E(\sigma), L(\sigma) \right]$. The minimum duration $D(\sigma)(t)$ and time warp use $TW(\sigma)(t)$ as a function of $t$ can be expressed as follows, where $D(\sigma)$ and $TW(\sigma)$ represent the minimum duration and time warp use, respectively:
    \begin{align}
        & D(\sigma)(t) = D(\sigma) + \max \{E(\sigma)-t), 0\} \label{eq:dispatch-early}\\
        & TW(\sigma)(t) = TW(\sigma) + \max \{t - L(\sigma), 0\} \label{eq:dispatch-late}
    \end{align}
\end{proposition}

The interval $\left[E(\sigma), L(\sigma)\right]$ is a time interval during which a vehicle should dispatch in order to minimize the duration and time warp use.
Equation~\eqref{eq:dispatch-early} describes how additional duration is incurred when dispatching earlier than $E(\sigma)$ (that is, due to additional waiting times).
Since we aim to minimize the total cost, we can always dispatch earlier without introducing additional costs.
Equation~\eqref{eq:dispatch-late} introduces an additional time warp penalty when starting later than $L(\sigma)$.
In particular, whenever a vehicle dispatches at its earliest dispatch time with $R^-(t) > L(\sigma)$, it is warped back in time to dispatch at time $L(\sigma)$.

We compute the total time warp involving penalties due to dispatch windows as follows:
\begin{align}
    TW_{DW}(\sigma) = TW(\sigma) +\max \{R^-(\sigma) - L(\sigma), 0\} + \max \{R^-(\sigma) - R^+(\sigma), 0\} \label{eq:time-warp-dispatch-windows}
\end{align}
The first term is the time warp penalty incurred from the concatenation scheme in Proposition~\ref{prop:concat-scheme}.
The second term is the time warp penalty incurred from starting later than $L(\sigma)$, which follows from Proposition~\ref{prop:starting-time}.
Finally, the third term adds an additional time warp penalty when the earliest dispatch time exceeds the dispatch time, i.e., whenever $R^-(\sigma) > R^+(\sigma)$.
Since there is a violation in the dispatch times, this expression introduces a ``forward'' time warp, meaning that the vehicle is warped from the latest dispatch time $R^+(\sigma)$ to the moment of earliest dispatch time $R^-(\sigma)$.
The total time warp~\eqref{eq:time-warp-dispatch-windows} is then used to evaluate local search moves.
For more details on the local search procedure, we refer to \cite{vidal2013hgs} and~\cite{wouda_pyvrp_2023}.

\subsection{Parameter tuning}
\label{sec:hgs-parameter-tuning}
This section outlines the parameter tuning procedure for the static solver and presents the resulting parameters.
PyVRP's HGS algorithm has a large number of parameters (see Table~\ref{tab:hgs-parameters}).
The default parameters of PyVRP are primarily optimized for solving larger VRP instances where more time is allowed.
Consequently, we use these parameters to compute the routing cost.
However, solving scenarios requires slightly different parameters to find good solutions within significantly shorter time limits of just a few seconds.

Because of the large number of parameters and operators in HGS, we opted for a prescriptive parameter tuning approach.
First, we derived reasonable lower and upper bounds for the values of each parameter.
Then, we split the parameters into five groups: (1) those related to population management, (2) those related to penalty management, (3) those related to neighborhood, (4) those related to the local search node operators, and finally, (5) those related to the local search route operators.
The genetic algorithm parameters were not tuned as they were irrelevant for solving scenario instances.
We generated 499 parameter configurations for a single group using a Latin hypercube design, and added the default setting as a control.
We used the set of 120 benchmark instances from Gehring and Homberger with 200 and 400 customers.
We solved each instance ten times using a different seed, for each of the 500 configurations.
We then took the best configuration and made that the default for the first group.
With this default, we then repeated this approach for the remaining groups.
Finally, we manually tuned the parameters using observations from the tuning process.
As shown in Table~\ref{tab:hgs-parameters}, the resulting parameter values for solving scenario instances are a smaller population size (allowing for faster convergence to good solutions), more frequent updates of the penalties, and the removal of many local search operators.
Tuning the parameters resulted in a respectable performance improvement when compared to the default parameters: on average, the solutions were 3\% better with a time limit of two seconds.

\begin{table}
  \caption{
  Parameters of PyVRP's solver, including the default parameter values and the parameter values for solving scenarios.
  Scenario parameter values are marked in bold if they differ from the default value.
}
  \label{tab:hgs-parameters}
\begin{tabular}{@{}lp{.57\linewidth}cc@{}}
\toprule
Category                           & Parameter                                         & Default & Scenario \\ \midrule
\multirow{2}{*}{Genetic algorithm} & Repair probability                                & 80\%  & 80\%  \\
                                   & Number of non-improving iterations before restart & 20000 & 20000 \\
\midrule
\multirow{6}{*}{Population}        & Minimum population size                           & 25    & \textbf{5}    \\
                                   & Population generation size                        & 40    & \textbf{3}    \\
                                   & Number of elite solutions                         & 4     & \textbf{2}     \\
                                   & Number of close solutions                         & 5     & \textbf{2}     \\
                                   & Lower bound diversity                             & 0.1   & 0.1   \\
                                   & Upper bound diversity                             & 0.5   & 0.5   \\
\midrule
\multirow{7}{*}{Penalty management}   & Initial capacity penalty                          & 20    & 20    \\
                                   & Initial time warp penalty                         & 6     & 6     \\
                                   & Repair booster                                    & 12    & 12    \\
                                   & Number of registrations between penalty updates   & 50    & \textbf{10}    \\
                                   & Penalty increase factor                           & 1.34  & \textbf{1.30}  \\
                                   & Penalty decrease factor                           & 0.32  & \textbf{0.50}  \\
                                   & Target feasible                                   & 0.43  & 0.43  \\
\midrule
\multirow{5}{*}{Neighborhood}     & Number of neighbors                                & 40    & 40    \\
                                   & Weight waiting time                               & 0.2     & 0.2   \\
                                   & Weight time warp                                  & 1.0     & 1.0   \\
                                   & Symmetric proximity                               & True  & True  \\
                                   & Symmetric neighbors                              & False  & False \\
\midrule
\multirow{10}{*}{Node operators}   & Include $(1, 0)$-exchange               & True  & True  \\
                                   & Include $(2, 0)$-exchange               & True  & \textbf{False}  \\
                                   & Include $(3, 0)$-exchange               & True  & \textbf{False}  \\
                                   & Include MoveTwoClientsReversed          & True  & \textbf{False}  \\
                                   & Include $(1, 1)$-exchange                   & True  & \textbf{False}  \\
                                   & Include $(2, 1)$-exchange                   & True  & \textbf{False}  \\
                                   & Include $(2, 2)$-exchange                   & True  & \textbf{False}  \\
                                   & Include $(3, 2)$-exchange                   & True  & \textbf{False}  \\
                                   & Include $(3, 3)$-exchange                   & True  & \textbf{False}  \\
                                   & Include 2-OPT operator                            & True  & True  \\
\midrule
\multirow{2}{*}{Route operators}  & Include RELOCATE* operator                        & True  & \textbf{False}  \\
                                   & Include SWAP* operator                            & True  & True  \\
\bottomrule
\end{tabular}
\end{table}

\clearpage
\section{Limited fleet}
\label{app:heterogeneous-fleet-extension}
In the problem formulation presented in Section~\ref{sec:problem-description}, we assume that there is an unlimited fleet of vehicles available at each epoch to serve all requests.
This assumption was introduced in the \textit{EURO meets NeurIPS} 2022 vehicle routing competition \citep{Kool2022a} to ensure that all requests could be served on time, but it clearly makes the problem less realistic since using more vehicles incurs no additional cost.
In this appendix, we discuss how to modify the original DDWP formulation to account for a limited number of cost-free vehicles, while ensuring that all requests can still be served on time.
Our approach is similar to \cite{vanHeeswijk2019}, making a distinction between a primary and secondary fleet of vehicles.

We assume that at each epoch $t$, there is a planned number of \textit{primary} vehicles that can be dispatched to serve requests.
Using vehicles from the primary fleet incurs no fixed costs, and unused primary vehicles from one epoch carry over to the next epoch.
In case the number of primary vehicles in an epoch is not enough to serve all requests, one can decide to hire vehicles from an external party (e.g., through crowd-sourcing).
We call vehicles from this fleet \textit{secondary}, and hiring additional secondary vehicles incur a fixed cost in addition to the routing cost.
In this modified problem, we can still ensure that all requests are served on time, but using more vehicles than planned is now penalized.
This change requires some adjustments to the static model and Markov decision process (MDP) since we must distinguish between a \textit{heterogeneous} fleet of vehicles.
The adjustments are described in the next subsections.

\subsection{Static model}
The modified static model extends the VRPTW-DW (\ref{subsec:static_vrptwdw}) as follows.
The fleet $\mathcal{K}$ now consists of the primary vehicles and the secondary vehicles, each with capacity $Q$.
Let $f_k \ge 0$ denote the fixed cost of using vehicle $k$, which is zero if $k$ is a primary vehicle and nonzero otherwise.
Each vehicle $k$ has an earliest time moment $E_k$ at which it is available for delivery.
This attribute is necessary to distinguish between vehicles that are available at the current epoch $t$, or some later epoch $t' > t$.
Indeed, whenever we sample a scenario instance for epoch $t' > t$, we must impose an earliest time moment $T_{t'}$ on primary vehicles that are only planned in future epochs to prevent them from serving requests in the current epoch.
The modified static model is then formulated as follows:
\begin{subequations}
\begin{align}
   \min_{x,\theta} \quad & \sum_{k \in \mathcal{K}} f_k \sum_{j \in \mathcal{N}} x_{0jk} + \sum_{k \in \mathcal{K}} \sum_{i, j \in \mathcal{V}} c_{ij} x_{ijk} \label{con:obj-hf} \\
    \text{s.t. } \quad & E_k \leq \theta_{0k} & \forall k \in \mathcal{K}, \label{con:vehicle-tw-early} \\
    &  (\ref{con:every-customer-assigned}) - (\ref{eq:constraint-1vrptw-rt}).
\end{align}
\label{milp:vrptw-dw-hf}
\end{subequations}

Objective~\eqref{con:obj-hf} minimizes the total vehicle fixed costs and the total routing cost.
Constraints~\eqref{con:vehicle-tw-early} ensure that a vehicle is dispatched no earlier than when it is available.
All constraints from the VRPTW-DW model~\eqref{milp:vrptw-dw} still apply.

\subsection{Markov decision process}
\subsubsection{State space}
The state is now a tuple $z_t = (s_t, K^p_t)$, where $s_t$ represents the set of requests that are known but not yet dispatched and $K^p_t$ denotes the number of primary vehicles at epoch $t$.

\subsubsection{Action space}
The action is now a tuple $b_t = (a_t, k^p_t)$, where $a_t \in \mathcal{A}(s_t)$ represents the subset of requests to dispatch and $k^p_t \in \{0, 1, \dots, K_t^p\}$ is the number of primary vehicles to use for dispatching the requests $a_t$.
The action space in state $z_t=(s_t, K^p_t)$ is thus defined as
\[\mathcal{B}(z_t) = \mathcal{A}(s_t) \times \{0, 1, \dots, K^p_t\} .\]

\subsubsection{State transition}
Let $P_t$ denote the primary vehicles that become available at epoch $t$.
The transition to the next state $z_{t+1}$ depends on the selected action $(a_t, k^p_t)$ and the unknown future realization $\omega_{t+1}$.
The state transitions from $z_t = (s_t, K^p_t)$ to $z_{t+1} = ((s_t \setminus a_t) \cup \omega_{t+1}, (K^p_t \setminus k^p_t) \cup P_{t+1})$.

\subsubsection{Optimality equation}
The direct cost $C(z_t, b_t)$ of a specific action $b_t=(a_t, k^p_t)$ in state $z_t$ is given by the cost of routing the requests in $a_t$ at time $T_t$ and using $k^p_t$ primary vehicles.
The direct cost is obtained by solving~\eqref{milp:vrptw-dw-hf} with $a_t$ as the set of requests and the fleet of vehicles $\mathcal{K}$ that consists of $k^p_t$ primary vehicles (with no fixed cost) and a large number of secondary vehicles (with fixed cost).
The objective of the DDWP is to select for each epoch $t \in \mathcal{T}$ a minimum cost action, such that the following (Bellman) optimality conditions are satisfied:
\begin{align}
    V(z_t) =
    \begin{dcases}
        \min_{b \in \mathcal{B}(z_t)} C(z_{t}, b) & \text{if } t = |\mathcal{T}|, \\
        \min_{b \in \mathcal{B}(z_t)} \big[ C(z_t, b)  + \mathbb{E}_{\omega_{t + 1}} [ V(z_{t+1}) ] \big] &  \text{otherwise}.
    \end{dcases}
\end{align}

\subsection{Numerical experiment}
We now evaluate the ICD methods on the modified DDWP with a limited primary fleet, following an identical experimental setup as used in Section~\ref{sec:experiments}.
In addition to what is described in the main body, we must determine the number of planned primary vehicles per epoch.
To this end, we first solve each instance assuming an unlimited primary fleet of vehicles using a greedy strategy that immediately dispatches all requests when revealed.
The resulting solution determines the number of planned primary vehicles for each epoch.
Indeed, it can be expected that operators plan a certain number of vehicles based on the expected demand, which is precisely what we aim to capture using this greedy strategy.
Furthermore, we assume that hiring a vehicle from the secondary fleet has a fixed upfront cost of two hours (recall that the routing cost equals the travel time).

We use the set of instances with time window characteristics and unimodal arrivals.
We solve these instances first assuming an unlimited fleet, and also assuming a limited fleet (obtained by the procedure outlined above).
We then compute the gap with respect to the hindsight solution obtained after 600 seconds, assuming an unlimited fleet.
Figure~\ref{fig:results-limited-fleet} plots the gap to the hindsight solution assuming an unlimited fleet and limited fleet for each method.
The results show that having only a limited number of primary vehicles based on a greedy strategy achieves a near-identical performance compared to the case with an unlimited fleet.
This means that limiting the vehicles (to a sufficient number to serve all requests) does not result in significantly more costs.

\begin{figure}
    \centering
    \includegraphics[width=0.8\linewidth]{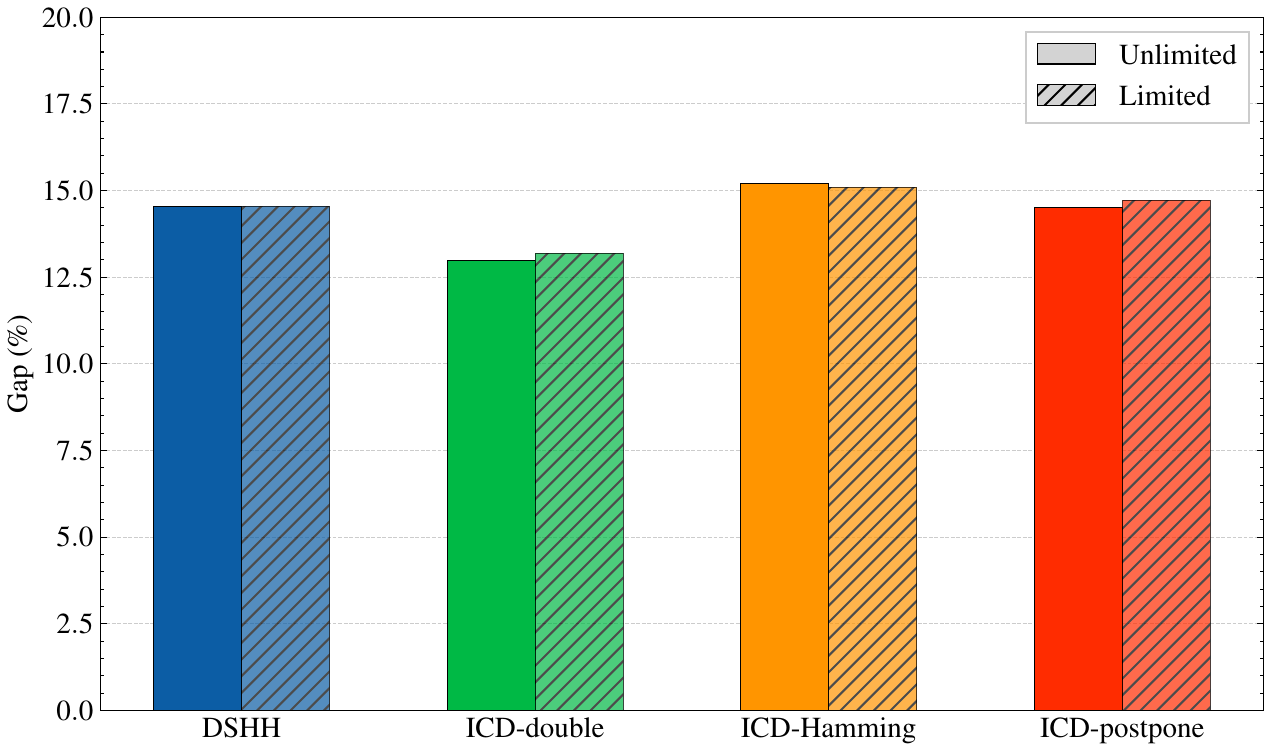}
    \caption{Average percentage gaps for different fleet configurations.}
    \label{fig:results-limited-fleet}
\end{figure}

Finally, it is not uncommon in the literature to use primary vehicles multiple times~\citep{Voccia2019,Klapp2018a,vanHeeswijk2019}.
This can be incorporated into the static problem formulation by redefining it as a multi-trip VRP, while also tracking the times at which the vehicle returns to the depot in the state description of the MDP.
We reserve this extension for future work.

\end{APPENDICES}

\clearpage
\bibliographystyle{informs2014trsc}
\bibliography{references} 

\end{document}